\def\Extended{extended}
\newif\ifextended
\ifx\Extended\Kaba\else\extendedtrue\fi
\newif\iftesting
\ifextended\testingtrue\fi
\ifextended
\documentclass[12pt]{article}\else
\documentclass{svproc}\fi
\ifextended
\renewcommand{\baselinestretch}{1.14}
\fi
\usepackage{url}

\ifextended
\usepackage[bookmarks=true,bookmarksnumbered=true,
  bookmarkstype=toc,
  pdftitle={Reflection principles, generic large cardinals, and the Continuum Problem},%
  pdfauthor={渕野 昌 (Sakaé Fuchino) and Andr\'e Ottenbreit Maschio Rodrigues},%
  pdfsubject={},%
  pdfkeywords={Continuum Problem, Laver-generically large cardinals, forcing 
  axioms, reflection principles}]{hyperref}
\hypersetup{
  colorlinks=true,
  linkcolor=blue,
  filecolor=blue,      
  citecolor=blue,
  urlcolor=cyan
}
\usepackage{graphicx, color}
\else
\usepackage{graphicx, color}
\fi
\usepackage{amsmath, amssymb, graphics}
\usepackage{bbm}% \mathbbm{} (for lowercase alhpabets) mathbb{} (for uppercase alphabets)
\usepackage{dsfont}
\usepackage{bbold}
\usepackage{mathtools}
\usepackage[normalem]{ulem}
\usepackage{setspace}
\def\memo#1{\iftesting\marginpar{{\normalsize\renewcommand{\baselinestretch}{0.4}\tiny%
			#1\par}}\else\fi}%
\newcount\minute	% Current minute within the hour
\newcount\hour		% Current hour (24-hour type)
\newcount\hourMins  % Temporary for taking hour modulo 60
\def\now%
{% Displays today's time as ``hh:mm''
  \minute=\time    % Number of minutes since midnight
  \hour=\time \divide \hour by 60 % Get hours
  \hourMins=\hour \multiply\hourMins by 60
  \advance\minute by -\hourMins % Hours modulo 60
  \zeroPadTwo{\the\hour}:\zeroPadTwo{\the\minute}%
}% --- \now ---
\def\zeroPadTwo#1%
{% Left zero pad of the argument to 2 digits.  The argument
  \ifnum #1<10 0\fi    % Conditionally output a zero
  #1%    Then output the argument
}% --- \zeroPadTwo ---
\newtheorem{Theorem}{Theorem}
\newtheorem{Lemma}[Theorem]{Lemma}
\newtheorem{Proposition}[Theorem]{Proposition}
\newtheorem{Corollary}[Theorem]{Corollary}
\ifextended
\newtheorem{claim}{Claim}[Theorem]
\newtheorem{problem}[Theorem]{Problem}
\fi
\iftesting
\newcommand{\Label}[1]{\label{#1}\marginpar{{\renewcommand{\baselinestretch}{0.6}\tiny 
		  #1}}}
\else
\newcommand{\Label}[1]{\label{#1}}
\fi%\iftesting
\newcounter{frml}[section]
\newcounter{frmla}[section]
\def\thefrml{{\arabic{section}.\arabic{frml}}}
\def\thefrmla{{$\aleph$\arabic{section}.\arabic{frmla}}}
\def\frmlabel#1{\refstepcounter{frml}{\def\baka{#1}\ifx\baka\empty\else\label{#1}\fi}%
{\rm({\thefrml})\hfill\hfill\hfill}}
\def\frmlabela#1{\refstepcounter{frml}{\def\baka{#1}\ifx\baka\empty\else\label{#1}\fi}%
{\rm({\thefrmla})\hfill\hfill\hfill}}
\iftesting
\def\xitem[#1]{\item[\frmlabel{#1}]%
	\marginpar{{\renewcommand{%
				\baselinestretch}{0.6}\tiny#1}}\ignorespaces}
\else
\def\xitem[#1]{\item[\frmlabel{#1}]\ignorespaces}
\fi
\def\xitemd[#1]#2{\item[(\ref{#1})$#2$\hfill\hfill\hfill]}
\def\xitema[#1]{\item[\frmlabela{#1}]\mbox{}%
	\iftesting\marginpar{{\renewcommand{%
				\baselinestretch}{0.6}\tiny#1}}\fi\ignorespaces}
\def\xitemsub[#1]#2{\item[\frmlabel{#1}$_{#2}$]\mbox{}%
	\iftesting\marginpar{{\renewcommand{%
				\baselinestretch}{0.6}\tiny#1}}\fi\ignorespaces}
\def\xxitem[#1][#2]{\item[(\ref{#1}{\makebox[1.4ex][c]{#2}})]\mbox{}%
	\iftesting\marginpar{{\renewcommand{%
				\baselinestretch}{0.6}\tiny\{#1\}\{#2\}}}\fi\ignorespaces}
\def\xitemof#1{{\rm({\ref{#1}})}}

\newenvironment{xitemize}{\begin{list}{}{\parsep=0.5\smallskipamount%
			\itemindent=-0.4ex%
			\itemsep=0.5\smallskipamount\itemindent=0em\leftmargin=3.8em\labelwidth=3em\labelsep=0.8em}}%
							 {\end{list}}
\def\assert#1{\noindent\makebox[4.8ex][r]{\rm(\makebox[2.2ex][c]{#1})}\ \ \ignorespaces}
\def\assertof#1{\makebox[3.5ex][c]{\rm(#1)}}%
\def\wassertof#1{\makebox[4.2ex][c]{\rm(\makebox[2.3ex][c]{#1})}}
\def\wassert#1{\noindent\makebox[4.2ex][c]{\rm(\makebox[2.3ex][c]{#1})}}
\newcommand{\prf}{{\bf Proof. }\ignorespaces}
\newcommand{\prfof}[1]{{\bf Proof of #1. }\ignorespaces}
\newcommand{\Qedof}[2]%
{\mbox{}\hspace*{\fill}{\ifextended$\square$\else\squareforqed\fi\,%
{\tiny\,\mbox{(#1~\number#2)}}}}
\newcommand{\sectionof}[1]{Section~\ref{#1}}
\newcommand{\Thmof}[1]{Theorem~\ref{#1}}
\newcommand{\Lemmaof}[1]{Lemma~\ref{#1}}
\newcommand{\Propof}[1]{Proposition~\ref{#1}}
\newcommand{\Corof}[1]{Corollary~\ref{#1}}
\newcommand{\Problemof}[1]{Problem~\ref{#1}}
\ifextended
\newcommand{\qed}{\mbox{}\hspace*{\fill}{$\square$}}
\fi
\newcommand{\qedofThm}{\Qedof{Theorem}{\theTheorem}}
\newcommand{\qedofLemma}{\Qedof{Lemma}{\theLemma}}
\newcommand{\qedofProp}{\Qedof{Proposition}{\theProposition}}
\newcommand{\qedofCor}{\Qedof{Corollary}{\theCorollary}}
\newcommand{\qedskip}{\medskip}
\newcommand{\ctenten}{,\mbox{}\hspace{0.08ex}{.}{.}{.}\hspace{0.1ex}}
\newcommand{\ctentenc}{,{}\linebreak[0]\hspace{0.04ex}{{.}{.}{.}\hspace{0.1ex},\,}\linebreak[0]}

\newcommand{\xmbox}[1]{ $\relax{\rm #1}\relax$ }
\newcommand{\utildeT}[1]{%
  \hbox to 0pt{\smash{$\mathop{\textstyle #1}\limits_{%
			\raisebox{0.4ex}[0pt]{$\scriptstyle\sim$}}$}\hss}%
  \relax\phantom{\mathord{{#1}_{\rule[-0.6ex]{0pt}{1pt}}}}}
\newcommand{\utildeS}[1]{%
	\hbox to 0pt{\smash{$\mathop{\scriptstyle #1}\limits_{%
				\raisebox{0.6ex}[0pt]{$\scriptscriptstyle\sim$}}$}\hss}%
	\relax\phantom{\mathord{{#1}_{\rule[-0.6ex]{0pt}{1pt}}}}}
\newcommand{\utildeSS}[1]{%
	\hbox to 0pt{$\mathop{\scriptscriptstyle #1}%
		\limits_{\scriptscriptstyle\sim}$\hss}%
		\relax\phantom{\underline{#1}}}
\newcommand{\utilde}[1]{%
	\mathchoice{\utildeT{#1}}{\utildeT{#1}}{\utildeS{#1}}{\utildeSS{#1}}}
\newcommand{\symb}[1]{{\mathord{\hspace{0.08em}\underbracket[0.6pt][2pt]{#1}}\hspace{0.08em}}}
\newcommand{\uta}{\utilde{a}}
\newcommand{\utb}{\utilde{b}}
\newcommand{\utu}{\utilde{u}}

\newcommand{\utpoP}{\utilde{\mathbb{P}}}
\newcommand{\utpoQ}{\utilde{\mathbb{Q}}}

\newcommand{\utgenH}{\utilde{\mathbb{H}}}
\newcommand{\utD}{\utilde{D}}
\newcommand{\utS}{\utilde{S}}
\newcommand{\bbone}{{\mathord{\mathbb{1}}}}

\newcommand{\circleq}{\mathrel{{\leqslant}%
		\hspace{-0.86ex}{\lower-0.53ex\hbox{$\scriptscriptstyle\circ$}}}}
\newcommand{\gmA}{\mathfrak{A}}
\newcommand{\gmB}{\mathfrak{B}}

\newcommand{\calA}{{\mathcal A}}

\newcommand{\calC}{{\mathcal C}}
\newcommand{\calD}{{\mathcal D}}
\newcommand{\calE}{{\mathcal E}}

\newcommand{\calG}{{\mathcal G}}
\newcommand{\calH}{{\mathcal H}}

\newcommand{\calL}{{\mathcal L}}

\newcommand{\calP}{{\mathcal P}}
\newcommand{\calQ}{{\mathcal Q}}

\newcommand{\calS}{{\mathcal S}}
\newcommand{\calT}{{\mathcal T}}

\newcommand{\st}{such that}
\newcommand{\wrt}{with respect to}
\newcommand{\ws}{winning strategy}
\newcommand{\Wolog}{Without loss of generality}

\newcommand{\Tfae}{The following are equivalent}
\newcommand{\tfae}{the following are equivalent}

\newcommand{\IMPLIES}{\ $\Rightarrow$\ }
\newcommand{\po}{po\-set}
\newcommand{\pos}{po\-sets}
\newcommand{\poP}{\mathbb{P}}
\newcommand{\poQ}{\mathbb{Q}}
\newcommand{\genG}{\mathbb{G}}
\newcommand{\genH}{\mathbb{H}}
\newcommand{\condp}{\mathbb{p}}
\newcommand{\condq}{\mathbb{q}}
\newcommand{\condr}{\mathbb{r}}

\newcommand{\utcondq}{\utilde{\mathbb{q}}}

\newcommand{\uniV}{{\sf V}}

\newcommand{\ZFC}{{\sf ZFC}}
\newcommand{\CH}{{\sf CH}}
\newcommand{\MA}{{\sf MA}}
\newcommand{\PFA}{{\sf PFA}}
\newcommand{\MM}{{\sf MM}}
\newcommand{\RP}{{\sf RP}}
\newcommand{\RC}{{\sf RC}}
\newcommand{\FRP}{{\sf FRP}}
\newcommand{\SDLS}{{\sf SDLS}}
\newcommand{\PKL}{{\sf PKL}}
\newcommand{\intnl}{int}
\newcommand{\GRP}{{\sf GRP}}
\newcommand{\LT}{{<}\,}
\newcommand{\LE}{{\leq}\,}
\newcommand{\varin}{\mathrel{\varepsilon}}
\newcommand{\notvarin}{\mathrel{\not\!\varepsilon\,}}
\newcommand{\setof}[2]{\{#1\,:\,#2\}}
\newcommand{\ssetof}[1]{\{#1\}}
\newcommand{\pairof}[1]{\langle#1\rangle}
\newcommand{\forces}[2]{\,\|\hspace{-.35ex}\mbox{\sf--}_{\,#1\,}%
\mbox{\rm``}\,#2\,\mbox{\rm''}}
\newcommand{\mapping}[3]{#1:#2\rightarrow #3}
\newcommand{\elembed}[3]{#1:#2\stackrel{\preccurlyeq\hspace{0.8ex}}{\rightarrow}#3}
\newcommand{\crit}{\mbox{\it crit\/}}
\newcommand{\fnsp}[2]{\mbox{}^{{#1}\hspace{-0.02em}}#2}
\newcommand{\imageof}{{}^{\,{\prime}{\prime}}}
\newcommand{\seqof}[2]{\langle#1\,:\,#2\rangle}
\newcommand{\modelof}[1]{\models\!\mbox{\rm``\,}#1\mbox{\rm''}}
\newcommand{\Pkl}[2]{\ifx\bakakaba#1\bakakaba\ifx\bakakaba#2\bakakaba{\mathcal 
    P}_\kappa(\lambda)\else{\mathcal P}_\kappa(#2)\fi\else{\mathcal P}_{#1}(#2)\fi}
\newcommand{\cardof}[1]{\mathopen{|\,}#1\mathclose{\,|}}
\newcommand{\psof}[1]{{\mathcal P}\/(#1)}
\newcommand{\concat}{\mathop{{}^{\frown}}}
\newcommand{\cf}{\mathop{\mbox{\it cf\/}}}
\newcommand{\otp}{\mathop{\mbox{\it otp\/}}}
\newcommand{\restr}{\restriction}
\newcommand{\Col}{{\rm Col}}
\newcommand{\continuum}{2^{\aleph_0}}
\begin{document}
\ifextended\else
\mainmatter              % start of a contribution
\fi
\title{Reflection principles, generic large cardinals, and the Continuum Problem}
\ifextended\else
\titlerunning{The Continuum Problem}  % abbreviated title (for running head)
\fi
\ifextended
\author{Saka\'e Fuchino${}^{(1)}$ \and Andr\'e Ottenbreit Maschio Rodrigues${}^{(2)}$ }
\else
\author{Saka\'e Fuchino\inst{1} \and Andr\'e Ottenbreit Maschio Rodrigues\inst{2}}
\fi
\ifextended\else
\authorrunning{Saka\'e Fuchino and Andr\'e Ottenbreit} % abbreviated author list (for running head)
\fi
\ifextended\else
\tocauthor{Saka\'e Fuchino and Andr\'e Ottenbreit Ottenbreit Maschio Rodrigues}
\fi
\ifextended
\else
\institute{Kobe University, Kobe, Japan\\
  \email{fuchino@diamond.kobe-u.ac.jp},\\
  WWW home page:\ \ifextended\href{https://fuchino.ddo.jp/}{\texttt{https://fuchino.ddo.jp/}}\else
  \texttt{https://fuchino.ddo.jp/}\fi
\and
Kobe University, Kobe, Japan\\
\email{andreomr@gmail.com}}
\fi
%% Universit\'{e} de Paris-Sud,
%% Laboratoire d'Analyse Num\'{e}rique, B\^{a}timent 425,\\
%% F-91405 Orsay Cedex, France}

\maketitle              % typeset the title of the contribution
\ifextended
\renewcommand{\thefootnote}{}
\footnotetext{${}^{(1)}$ Kobe University, Kobe, Japan\quad
  Email: {\tt fuchino@diamond.kobe-u.ac.jp}\\
  \hspace{1.1cm}Webpage: \href{https://fuchino.ddo.jp/}{\texttt{https://fuchino.ddo.jp/}}\smallskip\\
  \hspace{0.62cm}%
  ${}^{(2)}$ Kobe University, Kobe, Japan\quad
  Email: {\tt andreomr@gmail.com}  
}
\renewcommand{\thefootnote}{\arabic{footnote}}
\fi
\begin{abstract}
Strong reflection principles with the reflection cardinal $\leq\aleph_1$ or 
$\LT\continuum$ imply that the size of the continuum is either $\aleph_1$ or 
$\aleph_2$ or very large. Thus, the stipulation, that a strong reflection 
principle should hold, seems to support the trichotomy on the possible size of the continuum.

In this article, 
we examine the situation with the reflection principles and related notions of generic large 
cardinals.

\ifextended
\mbox{}\medskip

\else
\fi
\end{abstract}
\ifextended
%% contents 
\begin{quotation}
	\footnotesize
	\noindent
	\centerline{\normalsize\tt Contents\hspace{6em}\mbox{}}\mbox{}\\
       {\mbox{}\hspace{-1.6em}\tt\makebox[3.4ex][l]{\ref{multiverse}.}%
         G\"odel's Program and the set-theoretic multiverse}\ \ \dotfill\ \ \pageref{multiverse}\\ 
       {\mbox{}\hspace{-1.6em}\tt\makebox[3.4ex][l]{\ref{reflection}.}%
         Reflection Principles}\ \ \dotfill\ \ \pageref{reflection}\\ 
       {\mbox{}\hspace{-1.6em}\tt\makebox[3.4ex][l]{\ref{game-refl}.}%
         Game Reflection Principles and generically large cardinals}\ \ \dotfill\ \ \pageref{game-refl}\\ 
       {\mbox{}\hspace{-1.6em}\tt\makebox[3.4ex][l]{\ref{lt+le}.}%
         Simultaneous reflection down to $\LT\continuum$ and $\LE\continuum$}\ \ \dotfill\ \ \pageref{lt+le}\\ 
       {\mbox{}\hspace{-1.6em}\tt\makebox[3.4ex][l]{\ref{large-conti}.}%
         Reflection principles under large continuum}\ \ \dotfill\ \ \pageref{large-conti}\\ 
       {\mbox{}\hspace{-1.6em}\tt\makebox[3.4ex][l]{\ref{Laver}.}%
         Laver-generic large cardinals}\ \ \dotfill\ \ \pageref{Laver}\\ 
       {\mbox{}\hspace{-1.6em}\tt\makebox[3.4ex][l]{\ref{epilogue}.}%
         Some open problems and Takeuti's account on G\"odel's contribution to 
         the continuum problem}\ \ \dotfill\ \ \pageref{epilogue}\\ 
       %% {\mbox{}\hspace{-1.6em}\tt\makebox[3.4ex][l]{\ref{}.}%
       %%   \ \ \dotfill\ \ \pageref{}\\ 
       %% {\mbox{}\hspace{-1.6em}\tt\makebox[3.4ex][l]{\ref{}.}%
       %%   \ \ \dotfill\ \ \pageref{}\\ 
       %% {\mbox{}\hspace{-1.6em}\tt\makebox[3.4ex][l]{\ref{}.}%
       %%   \ \ \dotfill\ \ \pageref{}\\ 
       %% {\mbox{}\hspace{-1.6em}\tt\makebox[3.4ex][l]{\ref{}.}%
       %%   \ \ \dotfill\ \ \pageref{}\\ 
       {{\mbox{}\hspace{-1.6em}\tt\makebox[3.4ex][l]{}}%
         \tt
         References}\ \ \dotfill\ \ \pageref{ref}
\end{quotation}
\addtocounter{footnote}{1}
\footnotetext{The present text is an extended version of the paper with the same 
  title (last revised on: \today\ (\now\ JST)). Some details and proofs missing 
  in the version for publication, {\tt as well as the corrections done after the final 
  submitted version are 
  written in typewriter font}. The most up-to-date file of the present text with correct 
  typesetting of East Asian fonts is 
  downloadable as:
  \\\href{https://fuchino.ddo.jp/papers/refl\_principles\_gen\_large\_cardinals\_continuum\_problem-x.pdf}{%
    \scalebox{0.83}[1]{\tt 
    https://fuchino.ddo.jp/papers/refl\_principles\_gen\_large\_cardinals\_continuum\_problem-x.pdf}}\smallskip}
\else
\addtocounter{footnote}{1}
\footnotetext{An extended version of the manuscript of the paper with some more 
  details and proofs is downloadable as:
  \\\scalebox{0.83}[1]{\tt 
    https://fuchino.ddo.jp/papers/refl\_principles\_gen\_large\_cardinals\_continuum\_problem-x.pdf}}
\fi%ifextended
\addtocounter{footnote}{1}
\footnotetext{The authors would like to thank Hiroshi Sakai and the anonymous 
  referee for many valuable comments.}

\section{G\"odel's Program and large cardinals}
\Label{multiverse}
The Continuum Problem has been considered to be one of the central problems in set 
theory. Georg Cantor tried till the end of his mathematical carrier to prove 
his ``theorem'' which claims, formulated in present 
terminology, the continuum, the cardinality 
$2^{\aleph_0}$ of the set of all real numbers, is the first uncountable cardinal
$\aleph_1$. This statement is now called the Continuum Hypothesis (\CH). 
By G\"odel \cite{AC-GCH1}, \cite{AC-GCH2} \cite{V=L}, and Cohen \cite{cohen1}, 
\cite{cohen2}, \cite{cohen-book}, it is proven hat \CH\ is independent from the 
axiom system \ZFC\ of Zermelo-Fraenkel set theory with the Axiom of 
Choice.\footnote{Due to the Incompleteness Theorems, if we would like to 
  formulate this statement precisely, we have to put it under the assumption that 
  \ZFC\ is consistent (which we not only assume but do believe). } 

Although the majority of the non-set theorists apparently believes that the results by G\"odel 
and Cohen were the final solutions of the Continuum Problem, G\"odel maintained in 
\cite{what-is} that the conclusive solution to the problem is yet to be obtained in 
that a ``right'' extension of \ZFC\ will be found which will decide the size of the 
continuum. Today the research program of searching for possible legitimate extensions of 
\ZFC\ to settle the Continuum Problem is called G\"odel's 
Program. Now that, besides \CH,  a multitude of mathematically significant 
statements is known to be independent from \ZFC, the 
program should aim to decide not only the size of the continuum but also many of 
these independent mathematical statements. 
For modern views on G\"odel's Program, the reader may consult e.g.\ 
Bagaria \cite{bagaria},  
Steel \cite{steel}.

G\"odel suggested in \cite{what-is} that the large cardinal axioms are good 
candidates of axioms to be added to the axiom system \ZFC. Unfortunately large 
cardinals do not decide the size of the continuum which G\"odel also admits in 
the postscript to \cite{what-is} added in 1966. 
Nevertheless, it is known today that some notable structural aspects of the 
continuum like the Projective Determinacy are decided under the existence of 
certain large large cardinals.

In this paper, we discuss about a new notion of generic large cardinals introduced in 
Fuchino, Ottenbreit and Sakai \cite{II} and called there Laver-generic large 
cardinals (see \sectionof{Laver} below). 
Reasonable instances of 
(the existential statement of a) Laver-generic large cardinal decide the size of the continuum to be either 
$\aleph_1$ or $\aleph_2$ or fairly large. We show that these three 
possible scenarios of Laver-generic large cardinal are in accordance with 
respective strong reflection properties with reflection cardinal $\LT\aleph_2$ or $\LT\continuum$. 

In connection with the view-point of set-theoretical multiverse (see 
Fuchino \cite{multiverse}), our trichotomy theorems, or some further developments of them, 
have certain possibility to become the final answer to the Continuum Problem. As is 
well-known, Hugh Woodin is creating a theory which should 
support $\CH$ from the point of view of what should hold in a canonical model 
of the set theory. It should be emphasized that our trichotomy is not directly in 
contradiction with the possible outcome 
of his research program. In any case, it should be mathematical results in the 
future which should decide the matter definitively (if ever?). 

\ifextended
{\tt The authors learned the following ``bullet''-notation from Asaf Karagila at the RIMS set 
  theory workshop 2019. If $t(x_0\ctenten)$ is a term in some expansion of the 
  language of set theory by definitions then for a \po\ $\poP$ 
  and $\poP$-names $\uta_0$\ctenten, $t(\uta_0\ctenten)^\bullet$ denotes the 
  standard $\poP$-name $\utu$ \st\
  \begin{xitemize}
  \xitema[aa-0] 
    $\utu[\genG]=t^{\uniV[\genG]}(\uta_0[\genG]\ctenten)$. 
  \end{xitemize}
  For example
  $\pairof{\uta,\utb}^\bullet$ is denoted as ${\rm op}(\uta,\utb)$ in Kunen 
  \cite{kunen-book}. $t(\uta_0\ctenten)$ may have infinitely may parameters. For 
  example if $\uta_\xi$, $\xi<\delta$ is a sequence of $\poP$-names in the ground 
  model, $\setof{\uta_\xi}{\xi<\delta}^\bullet$ may be introduced as 
  the $\poP$-name $\setof{\pairof{\uta_\xi,\bbone_\poP}}{\xi<\delta}$. 
}\fi %%\ifextended
%%{\tt 

\section{Reflection Principles}
\Label{reflection}
The following type of mathematical reflection properties are considered in many 
different mathematical contexts. 
\begin{xitemize}
\xitem[refl-0]If a structure $\gmA$ in the class $\calC$ has the property $\calP$, then there is a 
  structure $\gmB$ in relation $\calQ$ to $\gmA$ \st\ $\gmB$ has 
  the cardinality $<\kappa$ and $\gmB$ also has the property $\calP$. 
\end{xitemize}
We shall call ``$<\kappa$'' in \xitemof{refl-0} above the 
{\it reflection cardinal} of the reflection property. If $\kappa$ is a successor 
cardinal $\mu^+$ we shall also say that the reflection cardinal is $\leq\mu$.

An example of an instance of \xitemof{refl-0} is, when $\calC=$ ``compact Hausdorff 
topological spaces'', $\calP=$ ``non-metrizable'', $\calQ=$ ``subspace'' and
$\kappa=\aleph_2$, that is, with the reflection cardinal $\leq\aleph_1$. In this 
case,  we obtain the statement: 
\begin{xitemize}
\xitem[refl-1] For any compact Hausdorff topological space, if $X$ is 
  non-metrizable, then there is a subspace $Y$ of $X$ of cardinality $<\aleph_2$ 
  \st\ $Y$ is also non-metrizable.
\end{xitemize}
This assertion is known to be a theorem in \ZFC\ (see Dow \cite{dow}).

If we extend the class $\calC$ in \xitemof{refl-1} to $\calC=$ ``locally compact Hausdorff 
space'', the statement thus obtained  
\begin{xitemize}
\xitem[refl-2] For any locally compact Hausdorff topological space, if $X$ is 
  non-metriz\-a\-ble, then there is a subspace $Y$ of $X$ of cardinality $<\aleph_2$ 
  \st\ $Y$ is also non-metrizable
\end{xitemize}
is no more a theorem in \ZFC: we can construct a counterexample to 
\xitemof{refl-2}, using a non-reflecting stationary subset $S$ of
$E^\kappa_\omega=\setof{\alpha<\kappa}{\cf(\alpha)=\omega}$ for some regular
$\kappa>\omega_1$ (Fuchino, Juh\'asz, Soukup, Szentmikl\'ossy and Usuba 
\cite{erice}). Note that $\square_\lambda$ for any uncountable 
$\lambda$ implies that there is such $S$ for $\kappa=\lambda^+$. 
In particular, \xitemof{refl-2} implies the total failure of the square principles 
and thus we need very large large cardinals to obtain the consistency of this 
reflection principle. Actually, a known consistency proof of this principle 
requires the existence of a strongly compact cardinal\footnote{The existence of a strongly compact cardinal 
is enough to force Rado's Conjecture discussed below and Rado's Conjecture 
implies the reflection statement \xitemof{refl-2}.}. 

\xitemof{refl-2} is equivalent to the stationarity reflection principle 
called Fodor-type Reflection Principle (\FRP) introduced in 
\cite{erice}.\footnote{Here, we are not only talking about 
  equiconsistency but really about equivalence over \ZFC.}
This principle can 
be formulated as follows (see \cite{more}). 

For a regular uncountable cardinal $\lambda$ and
$E\subseteq E^\lambda_\omega=\setof{\gamma\in\lambda}{\cf(\gamma)=\omega}$, a 
mapping  
$\mapping{g}{E}{[\lambda]^{\aleph_0}}$ is said to be a {\it ladder system on $E$} if, 
for all $\alpha\in E$, $g(\alpha)$ is a cofinal subset of $\alpha$ and
$\otp(g(\alpha))=\omega$.
\begin{xitemize}
\item[(\FRP)\,:] {\it For any regular $\lambda>\aleph_1$, stationary 
  $E\subseteq E^\lambda_\omega$, and a ladder system
  $\mapping{g}{E}{\lambda^{\aleph_0}}$ on $E$,
  there is an $\alpha^*\in E^\lambda_{\omega_1}$ \st\ 
  \begin{xitemize}
  \item[] $\setof{x\in[\alpha^*]^{\aleph_0}}{\sup(x)\in E,\,
    g(\sup(x))\subseteq x}$
  \end{xitemize}
  is stationary in $[\alpha^*]^{\aleph_0}$.
}
\end{xitemize}

Besides \xitemof{refl-2}, 
there are many mathematical reflection principles in the literature which have 
been previously studied rather separately but which are now all shown to be equivalent 
to \FRP\ and hence also equivalent to each other (see \cite{fuchino}, 
\cite{fuchino2}, \cite{fuchino-rinot}, \cite{more}). The equivalence of \xitemof{refl-2} 
to \FRP\ is established in \cite{more} via a further characterization of \FRP\ by 
non existence of a ladder system with a strong property of disjointness from 
which a counterexample to \xitemof{refl-2} (and other reflection properties 
proved to be equivalent to \FRP) can be constructed. Here 
we want to mention only one other reflection statement also equivalent to \FRP:

For a graph $G=\pairof{G,\calE}$, where $\calE\subseteq G^2$ is the adjacency 
relation of the graph, is said to be {\it of countable coloring number} if there is a 
well-ordering $\sqsubset$ on $G$ \st, for each $g\in G$,
$\setof{h\in G}{h\mathrel{\calE} g\mbox{ and }h\sqsubset g}$ is finite. 

The following assertion is also equivalent to \FRP\ (\cite{erice}, Fuchino, Sakai 
Soukup and Usuba \cite{more}):
\begin{xitemize}
\xitem[refl-3] For any graph $G$, if $G$ is not of countable coloring number, 
  then there is a subgraph $H$ of cardinality $<\aleph_2$ \st\ $H$ is neither of 
  countable coloring number. 
\end{xitemize}
In particular, it follows that the assertions \xitemof{refl-2} and 
\xitemof{refl-3} are equivalent to each other over \ZFC.

(Strong) Downward L\"owenheim Skolem Theorems of extended logics can be seen also as instances of the 
scheme \xitemof{refl-0}. The following is a theorem in \ZFC:

\begin{xitemize}
\item[$\SDLS(\calL(Q),\LT\aleph_2)$\,:] {\it For any uncountable first-order structure $\gmA$ in a countable 
  signature, there is an elementary submodel $\gmB$ of $\gmA$ \wrt\ the logic $\calL(Q)$ of 
  cardinality\footnote{The cardinality of a structure is defined to be the 
    cardinality of the underlying set.} $<\aleph_2$ where the quantifier $Q$ in a formula ``$Qx\,\varphi$'' is to be 
  interpreted as ``there are uncountably many $x$ \st\ $\varphi$''. }
\end{xitemize}

Adopting the notation of Fuchino, Ottenbreit and Sakai \cite{I}, let $\calL^{\aleph_0}_{stat}$ be the logic with 
monadic (weak) second order variable where the second order variables are to be 
interpreted as they are running over countable subsets of the structure in consideration. 
The logic has the built-in predicate $\varin$ where atomic formulas of the form
$x\varepsilon X$ is allowed for first and second order variables $x$ and $X$ 
respectively. The logic also has the unique second order quantifier $stat$ which is 
interpreted by
\begin{xitemize}
\xitem[refl-5] 
  for a structure
  $\gmA=\pairof{A\ctenten}$, $\gmA\models stat\,X\,\varphi[X\ctenten]$ holds if and only 
  if\\ $\setof{U\in[A]^{\aleph_0}}{\gmA\models\varphi[U\ctenten]}$ is stationary in
  $[A]^{\aleph_0}$. 
\end{xitemize}

Note that $\calL^{\aleph_0}_{stat}$ extends $\calL(Q)$ above, since $Qx\,\varphi$ 
can be expressed by \\$stat X\exists x\,(x\notvarin X\land\varphi)$. 

In $\calL^{\aleph_0}_{stat}$ we have two natural generalizations of the notion of 
elementary substructure. 
For (first order) structures $\gmA=\pairof{A\ctenten}$ and $\gmB=\pairof{B\ctenten}$ with $\gmB\subseteq\gmA$, let
\begin{xitemize}
\xitem[refl-6] $\gmB\prec_{\calL^{\aleph_0}_{stat}}\gmA$ if and only if, for all
  $\calL^{\aleph_0}_{stat}$-formula $\varphi=\varphi(x_0\ctentenc X_0\ctenten)$ 
  in the signature of $\gmA$, $b_0\ctenten\in B$, and 
  $U_0\ctenten\in [B]^{\aleph_0}$, we have\\
  $\gmB\models\varphi[b_0\ctentenc U_0\ctenten]$ $\Leftrightarrow$
  $\gmA\models\varphi[b_0\ctentenc U_0\ctenten]$. 
\xitem[refl-7] $\gmB\prec^-_{\calL^{\aleph_0}_{stat}}\gmA$ if and only if, for all
  $\calL^{\aleph_0}_{stat}$-formula $\varphi=\varphi(x_0\ctenten)$ in the signature of $\gmA$
  without any free second order variables, and $b_0\ctenten\in B$, 
  we have 
  $\gmB\models\varphi[b_0\ctenten]$ $\Leftrightarrow$
  $\gmA\models\varphi[b_0\ctenten]$. 
\end{xitemize}

By the remark after \xitemof{refl-5}, 
the following principles are generalizations of $\SDLS(\calL(Q),\LT\aleph_2)$:

\begin{xitemize}
\item[$\SDLS(\calL^{\aleph_0}_{stat},\LT\aleph_2)$\,:] {\it For any uncountable 
  first-order structure $\gmA$ in a countable signature, there is a submodel 
  $\gmB$ of $\gmA$ of cardinality $<\aleph_2$ \st\ $\gmB\prec_{\calL^{\aleph_0}_{stat}}\gmA$.}
  
\item[$\SDLS^-(\calL^{\aleph_0}_{stat},\LT\aleph_2)$\,:]  {\it For any uncountable 
  first-order structure $\gmA$ in a countable signature, there is a submodel 
  $\gmB$ of $\gmA$ of cardinality $<\aleph_2$ \st\ $\gmB\prec^-_{\calL^{\aleph_0}_{stat}}\gmA$.}
\end{xitemize}

M.\ Magidor noticed that $\SDLS^-(\calL^{\aleph_0}_{stat},\LT\aleph_2)$ implies 
\xitemof{refl-3} (see Magidor \cite{magidor2016}). By the equivalence of \xitemof{refl-3} 
to \FRP, we obtain
\begin{Theorem}
  \Label{T-refl-0}
  $\SDLS^-(\calL^{\aleph_0}_{stat},\LT\aleph_2)$ implies \FRP.\qed
\end{Theorem}

Actually, it is also easy to see that the stationarity reflection principle \RP\ 
(which is a strengthening of \RP\ 
in Jech \cite{millennium-book}) follows from
$\SDLS^-(\calL^{\aleph_0}_{stat},\LT\aleph_2)$.  

\FRP\ follows from our \RP\ (\cite{erice}) which is defined as follows:  

\begin{xitemize}
\item[\RP\,:] {\it For every regular $\lambda\geq\aleph_2$, stationary
  $S\subseteq[\lambda]^{\aleph_0}$, and $X\in[\lambda]^{\aleph_1}$, there is
  $Y\in[\lambda]^{\aleph_1}$ \st\ $\cf(Y)=\omega_1$, $X\subseteq Y$ and $S\cap[Y]^{\aleph_0}$ is 
  stationary in $[Y]^{\aleph_0}$.}
\end{xitemize}
Jech's \RP\ is just as our \RP\ as defined above but without demanding the 
property ``$\cf(Y)=\omega_1$'' for the reflection point $Y$.

\begin{Theorem}\Label{T-refl-1}
  $\SDLS^-(\calL^{\aleph_0}_{stat},\LT\aleph_2)$ implies \RP.
\end{Theorem}
{\bf Sketch of the proof.} Let $\lambda$, $S$, $X$ be as in the definition of 
\RP. Let $\mu>\lambda^{\aleph_0}$ be regular and
$\gmA=\pairof{\calH(\mu),\lambda,S, X, \in}$ where $\lambda$, $S$ 
and $X$ are thought to be interpretations of unary predicate symbols.  
Let $\gmB=\pairof{B\ctenten}$ be \st\ $B$ is of cardinality $\aleph_1$ and
$\gmB\prec^-_{\calL^{\aleph_0}_{stat}}\gmA$. Then $Y=\lambda\cap B$ is as 
desired. For example, $\cf(Y)=\omega_1$ follows from the fact that
$\gmB\models\psi$ by elementarity where $\psi$ is the $\calL^{\aleph_0}_{stat}$-sentence:
$stat\,X\,\exists y(y\mathbin{\symb{\in}}\symb{\lambda}\,\land\,\forall z\,
((z\varin X\land z\mathbin{\symb{\in}}\symb{\lambda})\rightarrow 
z\mathbin{\symb{\in}} y))$ 
where $\symb{\lambda}$ and $\symb{\in}$ are constant and  binary relation symbols 
corresponding to $\lambda$ and $\in$ in the structure $\gmA$. \qedofThm
\qedskip

By a theorem of Todor\v{c}evi\'c, \RP\ in the sense of Jech implies $\continuum\leq \aleph_2$ (see 
Theorem 37.18 in \cite{millennium-book}). Thus
\begin{Corollary}
  \Label{T-refl-2}
  $\SDLS^-(\calL^{\aleph_0}_{stat},\LT\aleph_2)$ 
  implies $2^{\aleph_0}\leq\aleph_2$. \qed
\end{Corollary}

In contrast to \Corof{T-refl-2}, \FRP\ does not put almost any restriction on the 
cardinality of the continuum since \FRP\ is preserved by ccc 
forcing (see \cite{erice}). 

A proof similar to that of \Thmof{T-refl-1} shows that 
$\SDLS^-(\calL^{\aleph_0}_{stat},\LT\aleph_2)$ implies the Diagonal Reflection 
Principle down to an internally club reflection point of cardinality $\LT\aleph_2$ 
of S.\ Cox \cite{cox}.  
Conversely, we can also easily prove that the Diagonal Reflection 
Principle down to an internally club reflection point of cardinality $\LT\aleph_2$ 
implies $\SDLS^-(\calL^{\aleph_0}_{stat},\LT\aleph_2)$. The internally clubness 
of the reflection point is used to guarantee that the internal interpretation of 
the stationary logic coincides with the external correct interpretation of the 
logic in the small substructure to make it an elementary substructure (in the 
sense of $\prec^-_{\calL^{\aleph_0}_{stat}}$)\vspace{-1\smallskipamount} of the original structure. Thus we 
obtain \assertof{1} of the following theorem. 
\begin{Theorem}{\rm (Theorem 1.1,\,\assertof{3} and \assertof{4} in \cite{I})}
  \Label{T-refl-3}\smallskip

  \wassert{1} $\SDLS^-(\calL^{\aleph_0}_{stat},\LT\aleph_2)$ 
  is equivalent to the Diagonal Reflection Principle down to an internally club 
  reflection point of cardinality $\LT\aleph_2$.\smallskip

  \wassert{2} $\SDLS(\calL^{\aleph_0}_{stat},\LT\aleph_2)$ is equivalent to
  $\SDLS^-(\calL^{\aleph_0}_{stat},\LT\aleph_2)$ plus \CH. \qed
\end{Theorem}

S.\,Cox proved in \cite{cox} that the Diagonal Reflection Principle down to an internally club 
reflection point of cardinality $\aleph_1$ follows from
$\MA^{+\omega_1}(\sigma\mbox{-closed})$. Thus, 

\begin{Corollary}\Label{C-refl-a}\wassertof{1} 
$\MA^{+\omega_1}(\sigma\mbox{-closed})$ implies $\SDLS^-(\calL^{\aleph_0}_{stat},\LT\aleph_2)$.\smallskip

\wassert{2}
$\MA^{+\omega_1}(\sigma\mbox{-closed})$ $+$ $\CH$ implies $\SDLS(\calL^{\aleph_0}_{stat},\LT\aleph_2)$ .\qed
\end{Corollary}

The reflection cardinal $\LT\aleph_2$ (or equivalently $\LE\aleph_1$) in the reflection 
principles above can be considered to be significant and even natural since, with this reflection 
cardinal, the reflection principles can be seen as  
statements claiming that the cardinality $\aleph_1$ is archetypical among 
uncountable cardinals, and hence that $\aleph_1$
already captures various phenomenon  
in uncountability in the sense that a certain type of properties of an uncountable 
structure can be reflected down to a substructure of the cardinality $\aleph_1$. 
From that point of view, it is interesting that one of the 
strongest reflection principles,  
namely the Strong Downward L\"owenheim-Skolem Theorem for stationary logic with this reflection cardinal 
implies \CH. 

In a similar way, we can also argue that the reflection 
with the reflection cardinal $\LT\continuum$ or $\LE\continuum$ should be 
regarded as significant and even natural since we can interpret the reflection with 
these reflection cardinals as 
a pronouncement of the richness of the continuum.

Let $\SDLS(\calL^{\aleph_0}_{stat},\LT\continuum)$ and
$\SDLS^-(\calL^{\aleph_0}_{stat},\LT\continuum)$  be the principles obtained from 
$\SDLS(\calL^{\aleph_0}_{stat},\LT\aleph_2)$ and $\SDLS^-(\calL^{\aleph_0}_{stat},\LT\aleph_2)$
by replacing ``$\LT\aleph_2$'' with ``$\LT\continuum$''.

\begin{Theorem}{\rm(Proposition 2.1, Corollary 2.3, Corollary 2.4 in \cite{II})}
  \Label{T-refl-4}\smallskip
  
  \wassert{1}
  $\SDLS^-(\calL^{\aleph_0}_{stat},\LT\continuum)$ 
  implies $\continuum=\aleph_2$. In particular, if $\continuum>\aleph_2$, then
  $\SDLS^-(\calL^{\aleph_0}_{stat},\LT\continuum)$ does not hold. 
\smallskip

  \wassert{2} $\SDLS(\calL^{\aleph_0}_{stat},\LT\continuum)$ is inconsistent. \qed

\end{Theorem}

Note that $\SDLS^-(\calL^{\aleph_0}_{stat},\LT\continuum)$ follows from
$\MA^{+\omega_1}(\sigma\mbox{-closed})$ $+$ $\neg\CH$ which is e.g.\ a 
consequence of $\PFA^{+\omega_1}$. 

Note that \Lemmaof{L-refl-0} implies that $\GRP^{\omega,\omega_1}(\LT\continuum)$ is also 
inconsistent. 

In contrast to the reflection down to $\LT\continuum$ whose strong version implies 
that the continuum is $\aleph_2$ (see \Thmof{T-refl-4},\,\assertof{2} above), the reflection down to $\LE\continuum$ does 
not exert any such restriction on the size of the continuum as we will see this 
in the next section. \memo{Lemma 4.11 in SDLS-x}

A slightly different type of reflection principle with reflection cardinal
$\LT\continuum$ implies that the continuum is very large. We will see this in 
\sectionof{large-conti}.

\section{Game Reflection Principles and generically large cardinals}
\Label{game-refl}
There is a further strengthening of 
$\SDLS^-(\calL^{\aleph_0}_{stat},\LT\aleph_2)$ which is called (Strong) Game 
Reflection Principle\footnote{In \cite{koenig}, B.\ K\"onig originally called the 
  principle introduced here the Strong Game  
  Reflection Principle and the local version of the principle the Game 
  Reflection Principle. } (\GRP) introduced in B.\ K\"onig \cite{koenig}. The 
following is a generalization of the principle:

For a regular uncountable cardinal $\mu$, 
a set $A$, and $\calA\subseteq\fnsp{\mu>}{A}$, 
$\calG^{\,\fnsp{\mu>}{A}}(\calA)$ is the following game of length $\mu$ for 
players I and I{}I. A match in $\calG^{\,\fnsp{\mu>}{A}}(\calA)$ looks like:\bigskip

\begin{tabular}{c|cccccccccc}
  I & \ \ $a_0$ && $a_1$ && $a_2$ && $\cdots$ &$a_\xi$ && $\cdots$\\
  \hline
  I{}I && $b_0$ && $b_1$ && $b_2$ & $\cdots$ && $b_\xi$ & $\cdots$
\end{tabular}\mbox{}\hfill ($\xi<\mu$)\qquad\mbox{}

\bigskip
\noindent
where 
$a_\xi$, $b_\xi\in A$ for $\xi<\mu$.

I{}I wins this match if
\begin{xitemize}
\xitem[] 
  $\seqof{a_\xi,b_\xi}{\xi<\eta}\in\calA$ and 
  $\seqof{a_\xi,b_\xi}{\xi<\eta}\concat\pairof{a_\eta}\not\in\calA$ for some
  $\eta<\mu$; or 
  $\seqof{a_\xi,b_\xi}{\xi<\mu}\in[\calA]$ 
\end{xitemize}
where
$\seqof{a_\xi,b_\xi}{\xi<\eta}$ denotes the sequence $f\in\fnsp{2\cdot\eta}{A}$ 
\st\ $f(2\cdot\xi)=a_\xi$ and $f(2\cdot\xi+1)=b_\xi$ for all $\xi<\eta$ and 
$[\calA]=\setof{f\in\fnsp{\mu}{A}}{f\restr\alpha\in\calA}$ for all $\alpha<\mu$. 

For regular cardinals $\mu$, $\kappa$ with $\mu<\kappa$
$\calC\subseteq[A]^{\LT\kappa}$ is said to be $\mu$-club if $\calC$ is cofinal in 
$[A]^{\LT\kappa}$ \wrt\ $\subseteq$ and closed \wrt\ the union of increasing 
$\subseteq$-chain of length $\nu$ for any regular $\mu\leq\nu<\kappa$. 

\begin{xitemize}
\item[$\GRP^{\LT\mu}(\LT\kappa)$:] {\it For any set $A$ of 
  regular cardinality $\geq\kappa$ and $\mu$-club
  $\calC\subseteq[A]^{\LT\kappa}$, if the player I{}I has no winning strategy in
  $\calG^{\,\fnsp{\mu>}{A}}(\calA)$ for some
  $\calA\subseteq\fnsp{\mu>}{A}$,  there is $B\in\calC$ \st\ the player I{}I has no 
  winning strategy 
  in $\calG^{\,\fnsp{\mu>}{B}}(\calA\cap\fnsp{\mu>}{B})$. }
\end{xitemize}

B.\ K\"onig's {\it Game Reflection Principle} (\GRP) is $\GRP^{<\omega_1}(<\aleph_2)$.

Sometimes, the following variation of the games and the principles is useful:
For a limit ordinal $\delta$, 
a set $A$, and $\calA\subseteq\fnsp{\delta\geq}{A}$, 
$\calG^{\,\fnsp{\delta\geq}{A}}(\calA)$ is the following game of length $\delta$ for 
players I and I{}I. A match in $\calG^{\,\fnsp{\delta\geq}{A}}(\calA)$ looks like:\bigskip

\begin{tabular}{c|cccccccccc}
  I & \ \ $a_0$ && $a_1$ && $a_2$ && $\cdots$ &$a_\xi$ && $\cdots$\\
  \hline
  I{}I && $b_0$ && $b_1$ && $b_2$ & $\cdots$ && $b_\xi$ & $\cdots$
\end{tabular}\mbox{}\hfill ($\xi<\delta$)\qquad\mbox{}

\bigskip
\noindent
where 
$a_\xi$, $b_\xi\in A$ for $\xi<\delta$.

I{}I wins this match if
\begin{xitemize}
\xitem[refl-8] 
  $\seqof{a_\xi,b_\xi}{\xi<\eta}\in\calA$ and 
  $\seqof{a_\xi,b_\xi}{\xi<\eta}\concat\pairof{a_\eta}\not\in\calA$ for some
  $\eta<\delta$; or 
  $\seqof{a_\xi,b_\xi}{\xi<\eta}\in\calA$ for all $\eta\leq\delta$.
\end{xitemize}
where
$\seqof{a_\xi,b_\xi}{\xi<\eta}$ is defined as above. 

For a limit ordinal $\delta$, and  uncountable regular cardinals $\mu$, $\kappa$ with $\delta\leq\mu<\kappa$, 
\begin{xitemize}
\item[$\GRP^{\delta,\mu}(\LT\kappa)$:] {\it For any set $A$ of 
  regular cardinality $\geq\kappa$ and $\mu$-club
  $\calC\subseteq[A]^{\LT\kappa}$, if the player I{}I has no winning strategy in
  $\calG^{\,\fnsp{\delta\geq}{A}}(\calA)$ for some
  $\calA\subseteq\fnsp{\delta\geq}{A}$,  there is $B\in\calC$ \st\ the player I{}I has no 
  winning strategy 
  in $\calG^{\,\fnsp{\delta\geq}{B}}(\calA\cap\fnsp{\delta\geq}{B})$. }
\end{xitemize}

The next Lemma follows immediately from the definitions: 

\begin{Lemma}
  \Label{T-refl-5} Suppose that $\delta$ and $\delta'$ are limit ordinals 
  and $\mu$, $\mu'$, $\kappa$, $\kappa'$ are regular cardinals \st\
  $\delta\leq\delta'<\mu\leq\mu'<\kappa$. Then we have

  \begin{xitemize}
  \xitem[refl-9] $\GRP^{<\mu'}(<\kappa)$ \IMPLIES\ $\GRP^{<\mu}(<\kappa)$ 
    \IMPLIES\ $\GRP^{\delta',\mu}(<\kappa)$ \IMPLIES\ $\GRP^{\delta,\mu}(<\kappa)$
    \ifextended\else\qed\fi
  \end{xitemize}
\end{Lemma}
\ifextended{\tt
\prf We prove the leftmost implication. The other implications 
can be proved similarly.

First, we assume that $\GRP^{\LT\mu'}(\LT\kappa)$ holds. Suppose 
$\calA\subseteq\fnsp{\mu>}{A}$ and $\calC\subseteq[A]^{\LT\kappa}$ is 
a $\mu$-club.
Let
\begin{xitemize}
  \xitema[refl-10]
  $\calA'=\setof{t\in\fnsp{\mu'>}{A}}{t\restr\alpha\in\calA\mbox{ for all }
  \alpha<\min\ssetof{\ell(t)+1,\mu}}$.
\end{xitemize}

Note that the player II has a \ws\ in $\calG^{\fnsp{\mu>}{A}}(\calA)$ if and only 
if the player II has a \ws\ in $\calG^{\fnsp{\mu'>}{A}}(\calA')$.

Suppose that the player II does not have any \ws\ in
$\calG^{\fnsp{\mu>}{A}}(\calA)$. Note that, since $\mu\leq\mu'$, $\calC$ is also
$\mu'$-club. By the equivalence above and by $\GRP^{\LT\mu'}(\LT\kappa)$, there 
is a $B\in\calC$ \st\ the player II does not have any \ws\ in
$\calG^{\fnsp{\mu'>}{B}}(\calA'\cap\fnsp{\mu'>}{B})$. By the definition 
of $\calA'$, it follows that the player II does not have any \ws\ in
$\calG^{\fnsp{\mu>}{B}}(\calA\cap\fnsp{\mu>}{B})$. This shows that we have
$\GRP^{<\mu}(<\kappa)$. \qedofLemma\qedskip
}\fi% \ifextended{\tt

\GRP\ is indeed a strengthening of $\SDLS(\calL^{\aleph_0}_{stat},\LT\aleph_2)$.
The following \Thmof{T-refl-6}, \Lemmaof{L-refl-0} and \Corof{C-refl-0} are 
slight generalizations of results in B.\ K\"onig \cite{koenig}.

\begin{Theorem}\Label{T-refl-6}{\rm (Theorem 4.7 in \cite{I})}
  Suppose that $\kappa$ is a regular uncountable cardinal \st\
  \begin{xitemize}
  \xitem[refl-11] $\mu^{\aleph_0}<\kappa$ for all $\mu<\kappa$, and
  \xitem[refl-12] $\GRP^{\omega,\omega_1}(<\kappa)$ holds.
  \end{xitemize}
  Then $\SDLS(\calL^{\aleph_0}_{stat},\LT\kappa)$ holds.\footnote{Actually we can 
    prove a slight strengthening of $\SDLS(\calL^{\aleph_0}_{stat},\LT\kappa)$ (see \cite{I}).}\qed
\end{Theorem}

\begin{Lemma}\Label{L-refl-0}{\rm (Lemma 4.2 in \cite{I})}
  For a regular cardinal $\kappa$, $\GRP^{\omega,\omega_1}(<\kappa)$ implies $2^{\aleph_0}<\kappa$. \qed
\end{Lemma}

Remember that \GRP\ is the principle $\GRP^{\LT\omega_1}(\LT\aleph_2)$. For a 
regular cardinal $\kappa>\aleph_1$ we shall write $\GRP(\LT\kappa)$ for
$\GRP^{\LT\omega_1}(\LT\kappa)$. Thus $\GRP$ is $\GRP(\LT\aleph_2)$. 

\begin{Corollary}\Label{C-refl-0}
  \wassertof{1} 
  \GRP\ implies $\SDLS(\calL^{\aleph_0}_{stat},\LT\aleph_2)$. %% In particular, 
  %% \GRP\ implies \CH.
  \smallskip

  \wassert{2} $\GRP(\LT(\continuum)^+)$ implies
  $\SDLS(\calL^{\aleph_0}_{stat},\LE\continuum)$. 
\end{Corollary}
\prf \assertof{1}: By \Lemmaof{L-refl-0}, \GRP\ implies \CH. Thus, under \GRP, 
\xitemof{refl-11} holds for $\kappa=\aleph_2$. By \Lemmaof{T-refl-5}, \GRP\ 
implies $\GRP^{\omega,\omega_1}(<\aleph_2)$. 
By \Thmof{T-refl-6}, it follows that $\SDLS(\calL^{\aleph_0}_{stat},\LT\aleph_2)$. 

\assertof{2}: Note that, for $\mu<(\continuum)^+$, 
$\mu^{\aleph_0}\leq\continuum<(\continuum)^+$ holds. By \Lemmaof{T-refl-5},
$\GRP(\LT(\continuum)^+)$ implies $\GRP^{\omega,\omega_1}(\LT(\continuum)^+)$. 
Thus, by \Thmof{T-refl-6}, it follows that
$\SDLS(\calL^{\aleph_0}_{stat},\LT(\continuum)^+)$, or 
$\SDLS(\calL^{\aleph_0}_{stat},\LE\continuum)$ in the other notation, holds. \\
\qedofCor
\qedskip

\GRP\ also implies another prominent reflection principle which is called Rado's 
Conjecture. 

We call a partial ordering $T=\pairof{T,\leq_T}$ a {\it tree} if the initial segment 
below any element is a well-ordering.  
A tree $T=\pairof{T,\leq_T}$ is said to be {\it special}\/ if it can be partitioned into
countably many antichains (i.e.\ pairwise incomparable sets). Note that every 
special tree has height $\leq\omega_1$. 

For a regular cardinal $\kappa>\aleph_1$, we define Rado's 
Conjecture with reflection cardinal $\LT\kappa$ as
\begin{xitemize}
\item[$\RC(\LT\kappa)$\,:] 
  {\it For any tree $T$, if $T$ is not special then there is $B\in[T]^{\LT\kappa}$ 
  \st\ $B$ (as the tree $\pairof{B,\leq_T\cap B^2}$) is not special.  }
\end{xitemize}

The original {\it Rado's Conjecture} (\RC) is $\RC(\LT\aleph_2)$.

\begin{Theorem}\Label{T-refl-7}{\rm(B.\ K\"onig \cite{koenig}, see also Theorem 
    4.3 in \cite{I})} For a regular cardinal
  $\kappa>\aleph_1$, $\GRP^{\LT\omega_1}(\LT\kappa)$ implies $\RC(\LT\kappa)$.\qed
\end{Theorem}

\FRP\ is also a consequence of \GRP. This is simply because \FRP\ follows from 
\RC\ (see \cite{rc-frp}). 

Game Reflection Principles are characterizations of certain instances of the 
existence of generically supercompact cardinals.

Let $\calP$ be a class of \pos. A cardinal $\kappa$ is said to be a {\it generically 
supercompact cardinal by $\calP$}, if, for any regular $\lambda$, there is a \po\ $\poP\in\calP$ 
\st, for any $(\uniV,\poP)$-generic filter $\genG$, there are classes $M$,
$j\subseteq\uniV[\genG]$ \st\ $M$ is an inner model of $\uniV[\genG]$, 
$\elembed{j}{\uniV}{M}$, $\crit(j)=\kappa$, 
$j(\kappa)>\lambda$ and $j\imageof\lambda\in M$.

\begin{Theorem}
  \Label{T-refl-8}{\rm(\cite{I})} For a regular uncountable $\kappa$, 
  \tfae:\smallskip

  \wassert{a} $2^{<\kappa}=\kappa$ and $\GRP^{\LT\kappa}(\LT\kappa^+)$ 
  holds.\smallskip

  \wassert{b} $\kappa^+$ is generically supercompact by $\LT\kappa$-closed \pos.\qed
\end{Theorem}
\begin{Corollary}\Label{C-refl-1}{\rm(B.\ K\"onig \cite{koenig})} \Tfae:\smallskip

  \wassert{a} \GRP\ holds.\smallskip

  \wassert{b} $\aleph_2$ is generically supercompact by $\sigma$-closed \pos. 
\end{Corollary}
\prf Assume that \GRP\ holds (remember that \GRP\ denotes
$\GRP^{\LT\omega_1}(\LT\aleph_2)$). Then, by \Corof{C-refl-0},\,\assertof{1}, 
$2^{\LT\aleph_1}=2^{\aleph_0}=\aleph_1$. Thus, by \Thmof{T-refl-8}, ``\assertof{a}
$\Rightarrow$ \assertof{b}'' for $\kappa=\aleph_1$, it follows that 
$\aleph_2=(\aleph_1)^+$ is generically supercompact by $\sigma$-closed forcing. 
The implication ``\assertof{b} $\Rightarrow$ \assertof{a}'' follows 
from ``\assertof{b} $\Rightarrow$ \assertof{a}'' of \Thmof{T-refl-8} for
$\kappa=\aleph_1$. 
\qedofCor

\section{Simultaneous reflection down to $\LT\continuum$ and $\LE\continuum$}
\Label{lt+le}
As we discussed in \sectionof{reflection}, the reflection down to $\LT\continuum$ as 
well as the reflection down to $\LE\continuum$ can be regarded as significant 
being principles which claim certain richness of the continuum.

One of the strong form of reflection principles with reflection cardinal
$\LT\continuum$ implies that the continuum is equal to $\aleph_2$ 
(\Thmof{T-refl-4},\,\assertof{2}) while there is a limitation on the possible types of 
reflection (\Thmof{T-refl-4},\,\assertof{3}).

In contrast, as we see below, the reflection down to $\LE\continuum$ can be 
established in one of its 
strongest forms without almost any restriction on the size of the continuum: \assertof{a} of \Thmof{T-refl-8} can be easily realized starting from a 
supercompact cardinal.

The following is well-known. 
\begin{Lemma}\Label{L-refl-1}{\rm(Lemma 4.10 in \cite{I})}
  If $\kappa$ is a supercompact and $\mu<\kappa$ is an uncountable regular 
  cardinal then for $\poP=\Col(\mu, \kappa)$ and $(\uniV,\poP)$-generic filter
  $\genG$, we have $\uniV[\genG]\modelof{\kappa=\mu^+\mbox{ and }\kappa\xmbox{ is generically 
  supercompact by }\LT\mu\xmbox{-closed \pos}}$.\qed
\end{Lemma}
%% \begin{Lemma}\Label{L-refl-2}{\rm(Lemma 4.12 in \cite{I})}
%%   Suppose that $\kappa^+$ is a generically measurable by
%%   $\LT\kappa$-closed forcing. Then $2^{\LT\kappa}=\kappa$ holds.\qed
%% \end{Lemma}

Suppose now that $\kappa_1$ is a supercompact cardinal and $\continuum$ is a 
regular cardinal. Let 
$\poQ=\Col(\continuum,\kappa_1)$ and let $\genH$ be a $(\uniV,\poQ)$-generic 
filter. By $\LT\continuum$-closedness of $\poQ$, we have
$(\continuum)^\uniV=(\continuum)^{\uniV[\genH]}$ and
$\uniV[\genH]\models\kappa_1=(\continuum)^+$. By \Lemmaof{L-refl-1},
$\uniV[\genH]\modelof{(\continuum)^+\xmbox{ is a generically supercompact cardinal by }\LT\continuum
  \xmbox{-closed \pos}}$. By \Thmof{T-refl-8}, it follows that
$\uniV[\genH]\modelof{\GRP^{\LT\continuum}(\LT(\continuum)^+)}$. 

By 
\Corof{C-refl-0},\,\assertof{2},
\Lemmaof{T-refl-5} and
\Thmof{T-refl-7}, we have, in particular, 
\begin{xitemize}
\xitem[refl-13] 
  $\uniV[\genH]\modelof{\SDLS(\calL^{\aleph_0}_{stat},\LE\continuum)\ \land\ 
    \RC(\leq\continuum)}$. 
\end{xitemize}

Note that the continuum can be forced to be practically anything of uncountable 
cofinality below $\kappa_1$ prior to the generic extension by $\poQ$.

The following \Propof{P-refl-0} should also belong to the folklore (for similar 
statements, see Theorem 4.1 in K\"onig and Yoshinobu \cite{KY} or Theorem 4.3 in 
Larson \cite{larson}). %% For a regular 
%% uncountable cardinal $\mu$, we shall call 
%% a \po\ $\poP$ {\it tightly $\LT\mu$-closed}\/ if any two compatible $\condp$,
%% $\condq\in\poP$ have their join $\inf\ssetof{\condp,\condq}$, and, 
%% for any decreasing sequence $\seqof{\condp_\xi}{\xi<\delta}$
%% of length $\delta<\mu$, there is the infimum of the sequence in $\poP$. Note 
%% that a L\'evy collapse with conditions of size $\LT\mu$ is tightly $\LT\mu$-closed. 

%% The following can be proved easily by induction on the size of the set $C$:

%% \begin{Lemma}
%%   \Label{L-refl-3}
%%   Suppose that $\poP$ is a tightly $\LT\mu$-closed \po\ 
%%   and $C\in[\poP]^{\LT\mu}$ is with the finite intersection property. Then there is $\condp\in\poP$ 
%%   \st\ $\condp=\inf(C)$.\qed
%% \end{Lemma}

Recall that, for a regular cardinal $\mu$, a \po\ $\poP$ is {\it $\LT\mu$-directed 
  closed}\/ if any downward directed subset of $\poP$  of cardinality $\LT\mu$ has a lower bound 
(in $\poP$). 

\begin{Proposition}
  \Label{P-refl-0} Suppose that $\MA^{+\omega_1}(\sigma\mbox{-closed})$ (or
  $\PFA^{+\omega_1}$, or $\MM^{+\omega_1}$, resp.) holds. If\/ $\poP$ is 
  $\LT\aleph_2$-directed closed, then we have
  \begin{xitemize}
  \xitem[refl-14] 
    $\forces{\poP}{\MA^{+\omega_1}(\sigma\mbox{-closed})\ \ 
    \mbox{(or }\PFA^{+\omega_1}\mbox{, or }\MM^{+\omega_1}\mbox{, resp.)}}$. 
  \end{xitemize}
\end{Proposition}
\prf We prove the case of $\MA^{+\omega_1}(\sigma\mbox{-closed})$. Other cases 
can be proved by the same argument. 

Suppose that $\poP$ is a $\LT\aleph_2$-directed closed \po\ and 
let $\utpoQ$, $\seqof{\utD_\alpha}{\alpha<\omega_1}$,
$\seqof{\utS_\beta}{\beta<\omega_1}$ be $\poP$-names \st\ 
\begin{xitemize}
\xitem[refl-15] $\forces{\poP}{%
  \begin{array}[t]{@{}l}
    \utpoQ\mbox{ is a }\sigma\mbox{-closed \po,}\\[\jot]
    \utD_\alpha\ (\alpha<\omega_1)\mbox{ is a dense subset 
      of }\utpoQ\mbox{ for all }\alpha<\omega_1,\mbox{ and}\\[\jot]
    \utS_\beta\ (\beta<\omega_1)\mbox{ is a }\poQ\mbox{-name of a stationary subset of }\omega_1\\[\jot]
    \mbox{for all }\beta<\omega_1}
  \end{array}$
\end{xitemize}

Let $\poP^*=\poP\ast\utpoQ$. For $\alpha<\omega_1$, let
\begin{xitemize}
\xitem[refl-16] 
  $D^*_\alpha=\setof{\pairof{\condp,\utcondq}\in\poP^*}{\condp
  \forces{\poP}{\utcondq\varin\utD_\alpha}}$.  
\end{xitemize}
For $\beta<\omega_1$, let
\begin{xitemize}
\xitem[refl-17] 
  $\utS^*_\beta=\setof{\pairof{\pairof{\condp,\utcondq},\check{\alpha}}}{
  \pairof{\condp,\utcondq}\in\poP^*,\,
  \condp\forces{\poP}{\utcondq\forces{\utpoQ}{\check{\alpha}\varin\utS_\beta}}}$.
\end{xitemize}
By the definition of $\poP^*$, $\seqof{D^*_\alpha}{\alpha<\omega_1}$, and
$\seqof{\utS^*_\beta}{\beta<\omega_1}$, 
the following is easy to show:
\begin{claim}
  $\poP^*$ is a $\sigma$-closed \po, $D^*_\alpha$ is a dense subset of $\poP^*$ 
  for all $\alpha<\omega_1$, and   $\utS^*_\beta$ is a $\poP^*$-name with
  $\forces{\poP^*}{\utS^*_\beta\xmbox{ is a stationary subset of }\omega_1}$ for all
  $\beta<\omega_1$. \\[-1.5\medskipamount]\qed\vspace{-1\smallskipamount}
\end{claim}

Let $\calD^*=\setof{D^*_\alpha}{\alpha<\omega_1}$. By
$\MA^{+\omega_1}(\sigma\mbox{-closed})$, there is a $\calD^*$-generic filter 
$\genG^*$ on $\poP^*$ \st\ $\utS^*_\beta[\genG^*]$ is a stationary subset of 
$\omega_1$ for all $\beta<\omega_1$. 

Let $\theta$ be a sufficiently large regular cardinal and let 
$M\prec\calH(\theta)$ be of cardinality $\aleph_1$ \st\ $\omega_1\subseteq M$ and
$M$ contains everything relevant (in particular, $\genG^*$, $D^*_\alpha$,
$\utS^*_\beta\in M$ for $\alpha$, $\beta<\omega_1$).

Let $\genG_0=\genG^*\cap M$ and let $\genG$ be the filter on $\poP^*$ generated by
$\genG_0$. By the choice of $M$, we have
$\utS^*_\beta[\genG^*]=\utS^*_\beta[\genG_0]=\utS^*_\beta[\genG]$. \vspace{0.6\smallskipamount}

Let
$G=\setof{\condp\in\poP}{\pairof{\condp,\utcondq}\in\genG\mbox{ for some }\utcondq}$. \vspace{0.6\smallskipamount}
Since $\cardof{G}\leq\cardof{M}<\aleph_2$ and $G$ is \smallskip downward 
directed, there is a lower bound $\condp_0\in\poP$ of $G$. 

Let
\begin{xitemize}
\xitem[refl-18] 
  $\utgenH=\setof{\pairof{\utcondq,\bbone_\poP}}{\pairof{\condp,\utcondq}\in\genG\mbox{ for some }\condp\in\poP}$.
\end{xitemize}
Then $\utgenH$ is a $\poP$-name and we have 
\begin{xitemize}
\xitem[refl-19] 
  $\condp_0\forces{\poP}{{}
  \begin{array}[t]{@{}l}
    \utgenH\mbox{ is a }\setof{\utD_\alpha}{\alpha<\omega_1}^\bullet
    \mbox{-generic filter on }\utpoQ\mbox{ \st}\\[\jot]
    \utS_\beta[\utgenH]\mbox{ is a stationary subset of }\omega_1
    \mbox{ for all }\beta<\omega_1}.\end{array}$
\end{xitemize}

Since the argument above can be also performed in $\poP\restr\condr$ instead of 
in $\poP$ for any $\condr\in\poP$. It follows that 
\begin{xitemize}
\xitem[refl-19-0] 
  $\forces{\poP}{{}
  \begin{array}[t]{@{}l}
    \mbox{there is a }\setof{\utD_\alpha}{\alpha<\omega_1}^\bullet
    \mbox{-generic filter }H\mbox{ on }\utpoQ\mbox{ \st}\\[\jot]
    \utS_\beta[H]\mbox{ is a stationary subset of }\omega_1
    \mbox{ for all }\beta<\omega_1}.\end{array}$\\
  \qedofProp
\end{xitemize}

\begin{Theorem}\Label{T-refl-9}
  Suppose that $\kappa$ and $\kappa_1$ with $\kappa<\kappa_1$ are two 
  supercompact cardinals. Then there is a generic 
  extension $\uniV[\genG\ast\genH]$ \st\
  \begin{xitemize}
  \item[] 
    $\uniV[\genG\ast\genH]\models\MM^{+\omega_1}$ $+$
    $\GRP^{\LT\continuum}(\leq\continuum)$.   
  \end{xitemize}
\end{Theorem}

Note that, by \Corof{C-refl-a},\,\assertof{1}, we have 
\begin{xitemize}
\item[] 
  $\uniV[\genG\ast\genH]\models\SDLS^-(\calL^{\aleph_0}_{stat},\LT\continuum)$ $+$ 
  $\GRP^{\LT\continuum}(\leq\continuum)$.
\end{xitemize}

\noindent
\prfof{\Thmof{T-refl-9}} Let $\uniV[\genG]$ be a standard model of \MM\ obtained 
by a reverse countable iteration of length $\kappa$ along with a fixed 
Laver-function $\kappa\,\rightarrow\,V_\kappa$. It is easy to see that $\uniV[\genG]$ also satisfies
$\MM^{+\omega_1}$. Note that we have $\uniV[\genG]\models\kappa=\aleph_2=\continuum$. 
In $\uniV[\genG]$, $\kappa_1$ is still supercompact. Thus, 
working in $\uniV[\genG]$, let $\poQ=\Col(\continuum,\kappa_1)$. 
Let $\genH$ be a $(\uniV[\genG], \poQ)$-generic filter. By \Propof{P-refl-0}, we 
have $\uniV[\genG\ast\genH]\models \MM^{+\omega_1}$. By \Lemmaof{L-refl-1} and 
\Thmof{T-refl-8}, we have
$\uniV[\genG\ast\genH]=(\uniV[\genG])[\genH]\models\GRP^{\LT\continuum}(\LE\continuum)$. 
\qedofThm

\section{Reflection principles under large continuum}
\Label{large-conti}
The continuum can be ``very large'' as a cardinal number. For example, this is 
the case in the model $\uniV[\genG]$ obtained by starting from a supercompact $\kappa$ and then 
adding $\kappa$ many Cohen reals. In this 
model, we have $\continuum=\kappa$ and there is a countably saturated normal fine filter over $\Pkl{}{}$ 
for all regular $\lambda\geq\kappa$. The last property of $\uniV[\genG]$ implies 
that $\kappa$ there is still fairly large 
(e.g. $\kappa$-weakly Mahlo and more, see e.g. Proposition 16.8 in Kanamori \cite{kanamori}).

If the ground model satisfies \FRP\ then 
$\uniV[\genG]$ also satisfies \FRP\ since \FRP\ is preserved by ccc extensions (see \cite{erice}). 
On the other hand, as we already have seen,
$\SDLS^-(\calL^{\aleph_0}_{stat},\LT\aleph_2)$ or even
$\SDLS^-(\calL^{\aleph_0}_{stat},\LT\continuum)$ is incompatible with large 
continuum. In particular, these reflection principles do not hold in our 
model $\uniV[\genG]$. 

A weakening of $\SDLS^-(\calL^{\aleph_0}_{stat},\LT\continuum)$ is compatible 
with large continuum. Let us begin with the diagonal reflection principle which 
characterizes the version of the strong downward L\"owenheim-Skolem theorem 
with reflection points of cardinality $<$ large continuum. 
The following is a weakening of Cox's Diagonal Reflection 
Principle down to an internally club reflection point. 

For regular cardinals $\kappa$, $\lambda$ with
$\kappa\leq\lambda$, let 
\begin{xitemize}
\item[$(\ast)^{\intnl+}_{\LT\kappa,\lambda}$:] {\it For any countable expansion 
$\tilde{\gmA}$ of $\pairof{\calH(\lambda),\in}$ and sequence 
$\seqof{S_a}{a\in\calH(\lambda)}$ \st\ $S_a$ is a stationary subset of 
$[\calH(\lambda)]^{\aleph_0}$ for all $a\in\calH(\lambda)$, there are 
  stationarily many 
  $M\in[\calH(\lambda)]^{\LT\kappa}$ \st\ 
  \begin{xitemize}
  \item[\wassertof{1}] $\tilde{\gmA}\restr M\prec\tilde{\gmA}$; and 
  \item[\wassertof{2}] $S_a\cap M$ is stationary in $[M]^{\aleph_0}$ 
    for all $a\in M$.
  \end{xitemize}}
\end{xitemize}
Note that \assertof{1} implies that $c\subseteq M$ holds for all
$c\in [M]^{\aleph_0}\cap M$. 

In the notation above, ``\intnl'' (internal) refers to the condition \assertof{2} 
in which not $S_a\cap [M]^{\aleph_0}$ but $S_a\cap M$ is declared to be 
stationary in $[M]^{\aleph_0}$; ``$+$'' refers to the condition that
$M\in[\calH(\lambda)]^{<\kappa}$ with \assertof{1} and \assertof{2} not only exists but there 
are stationarily many such $M$. 

That $(\ast)^{\intnl+}_{\LT\kappa,\lambda}$ is compatible 
with $\kappa=\continuum$ and it is arbitrarily large is seen in the following 
\Thmof{T-refl-10} together with \Lemmaof{L-PKL-0} below: 

\begin{Theorem}\Label{T-refl-10}{\rm (Theorem 2.10 in \cite{II})}% P-lt-conti-3
  Suppose that $\kappa$ is a generically supercompact 
  cardinal by
  %% $\kappa$-cc
  proper \pos. Then $(\ast)^{\intnl+}_{\LT\kappa,\lambda}$ holds for all 
  regular $\lambda\geq\kappa$. \qed
\end{Theorem}

Similarly to \Lemmaof{L-refl-1}, starting from a supercompact cardinal, it is easy to force that the continuum is 
generically supercompact cardinal by ccc-\pos. Let us call a \po\ $\poP$ 
appropriate for $\kappa$, if we have $j\imageof\poP\circleq j(\poP)$ for all supercompact 
embedding $j$ for $\kappa$. 
\begin{Lemma}\Label{L-PKL-0}
  If $\kappa$ is a supercompact and $\mu<\kappa$ is an uncountable regular 
  cardinal then for any $\LT\mu$-cc \po\ $\poP$ appropriate for $\kappa$, adding $\geq\kappa$ many 
  reals,  we have $\uniV[\genG]\modelof{\kappa\leq\continuum\mbox{ and }\kappa\xmbox{ is generically 
  supercompact by }\LT\mu\xmbox{-cc \pos}}$.\qed
\end{Lemma}

``$(\ast)^{\intnl+}_{\LT\kappa,\lambda}$ holds for all 
  regular $\lambda\geq\kappa$'' characterizes the strong downward 
L\"owenheim-Skolem theorem for internal interpretation of stationary logic 
defined in the following.

For a structure $\gmA=\pairof{A\ctenten}$ of a countable signature, an
$\calL^{\aleph_0}_{stat}$-formula
$\varphi=\varphi(x_0\ctentenc X_0\ctenten)$\footnote{As before, when we write
  $\varphi=\varphi(x_0\ctentenc X_0\ctenten)$, we always assume that the list $x_0\ctenten$ 
  contains all the free first order variables of $\varphi$ and $X_0\ctenten$ 
  all the free weak second order variables of $\varphi$.} and
$a_0\ctenten\in A$, $U_0\ctenten\in[A]^{\aleph_0}\cap A$, we define the internal 
interpretation of $\varphi(a_0\ctentenc U_0\ctenten)$ in $\gmA$ (notation:
$\gmA\models^{\intnl}\varphi(a_0\ctentenc U_0\ctenten)$ for
``$\varphi(a_0\ctentenc U_0\ctenten)$ holds internally in $\gmA$'') by induction 
on the construction of $\varphi$ as follows:

If $\varphi$ is ``$x_i\varin X_j$'' then
\begin{xitemize}
\xitem[internal-0] 
  $\gmA\models^{\intnl}\varphi(a_0\ctentenc U_0\ctenten)$\ \ $\Leftrightarrow$\ \
  $a_i\in U_j$
\end{xitemize}
for a structure $\gmA=\pairof{A\ctenten}$, $a_0\ctenten\in A$ and
$U_0\ctenten\in [A]^{\aleph_0}\cap A$. 

For first-order connectives and quantifiers in
$\calL^{\aleph_0}_{stat}$, the semantics ``$\models^{\intnl}$'' is defined exactly as for 
the first order ``$\models$''.

For 
an $\calL^{\aleph_0}_{stat}$ formula $\varphi$ with
$\varphi=\varphi(x_0\ctentenc X_0\ctentenc X)$, assuming that the notion of 
$\gmA\models^{\intnl}\varphi(a_0\ctentenc U_0\ctentenc U)$ has 
been defined for all $a_0\ctenten\in A$, $U_0\ctentenc U\in[A]^{\aleph_0}\cap A$, 
we stipulate 
\begin{xitemize}
\xitem[internal-1] 
  $\gmA\models^{\intnl}stat\,X\,\varphi(a_0\ctentenc U_0\ctentenc X)$\ \
  $\Leftrightarrow$\\
  \qquad$\setof{U\in[A]^{\aleph_0}\cap A}{\gmA\models^{\intnl}\varphi(a_0\ctentenc U_0\ctentenc U)}$
  is stationary in $[A]^{\aleph_0}$
\end{xitemize}
for a structure $\gmA=\pairof{A\ctenten}$ of a relevant signature,  
$a_0\ctenten\in A$ and $U_0\ctenten\in[A]^{\aleph_0}\cap A$.

For structures $\gmA$, $\gmB$ of the same signature 
with $\gmB=\pairof{B\ctenten}$ and $\gmB\subseteq\gmA$, 
we define
\begin{xitemize}
\xitem[internal-2] 
  {\it 
  $\gmB\prec^{\intnl}_{\calL^{\aleph_0}_{stat}}\gmA$\ \ $\Leftrightarrow$\\[2\jot]
  \qquad $\gmB\models^{\intnl}\varphi(b_0\ctentenc U_0\ctenten)$ if and only if 
  $\gmA\models^{\intnl}\varphi(b_0\ctentenc U_0\ctenten)$\\
  \qquad for all $\calL^{\aleph_0}_{stat}$-formulas $\varphi$ in the signature 
  of the structures with\\
  \qquad $\varphi=\varphi(x_0\ctentenc X_0\ctenten)$, $b_0\ctenten\in B$ and $U_0\ctenten\in[B]^{\aleph_0}\cap B$.}
\end{xitemize}

Finally, for a regular $\kappa>\aleph_1$, the internal strong downward 
L\"owenheim-Skolem Theorem 
$\SDLS^{\intnl}_+(\calL^{\aleph_0}_{stat},\LT\kappa)$ 
is defined by 
\begin{xitemize}
\item[$\SDLS^{\intnl}_+(\calL^{\aleph_0}_{stat},\LT\kappa)$:]  {\it For any structure $\gmA=\pairof{A\ctenten}$ of 
  countable signature with $\cardof{A}\geq\kappa$, there are {\it stationarily many} 
  $M\in[A]^{\LT\kappa}$ \st\\ $\gmA\restr M\prec^{\intnl}_{\calL^{\aleph_0}_{stat}}\gmA$. }
\end{xitemize}

Similarly to the $+$ in ``$(\ast)^{\intnl+}_{\LT\kappa,\lambda}$'',  
`$+$' in ``$\SDLS^{\intnl}_+(\calL^{\aleph_0}_{stat},\LT\kappa)$'' refers to 
the existence of ``stationarily many'' reflection points $M$.
This additional condition can be dropped if
$\kappa=\aleph_2$. This is because the quantifier $Qx\,\varphi$ defined by
$stat\,X\exists x\,(x\notvarin X\land\varphi$, 
$\gmA\models^{\intnl}Qx\,\varphi(x\ctenten))$ still implies that ``there are 
uncountably many $a\in A$ with $\varphi(a\ctenten)$''. 
Note that, if $\gmA\models^{\intnl}\neg stat\,X\,(x\equiv x)$, for a structure
$\gmA=\pairof{A\ctenten}$,  we can easily find even club many  
$X\in[A]^{<\kappa}$ for any regular $\aleph_1\leq\kappa\leq\cardof{A}$ \st\ 
$\gmA\restr X\prec^{\intnl}_{\calL^{\aleph_0}_{stat}}\gmA$. 

\begin{Proposition}\Label{P-internal-0}{\rm (Proposition 3.1 in \cite{II})}
  For a regular cardinal $\kappa>\aleph_1$, \tfae:\smallskip

  \assert{a} $(\ast)^{\intnl+}_{\LT\kappa,\lambda}$ holds for all 
  regular $\lambda\geq\kappa$.\smallskip

  \assert{b} $\SDLS^{\intnl}_+(\calL^{\aleph_0}_{stat},\LT\kappa)$ holds.\qed
\end{Proposition}

Although $\SDLS^{\intnl}_+(\calL^{\aleph_0}_{stat},\LT\continuum)$ is compatible with 
large continuum, as a weakening 
of $\SDLS^-(\calL^{\aleph_0}_{stat},\LT\continuum)$, this principle does not 
imply the largeness of the continuum. The strong L\"owenheim-Skolem theorem for 
the following variation of stationary logic does. 

For sets $s$ and $t$ we denote with $\Pkl{s}{t}$ the set
$[t]^{\LT\cardof{s}}=\setof{a\in\psof{t}}{\cardof{a}<\cardof{s}}$. We say
$S\subseteq\Pkl{s}{t}$ is stationary if it is stationary in the sense of Jech 
\cite{millennium-book}. 

The logic $\calL^{\PKL}_{stat}$ has a built-in unary predicate symbol
$\symb{K}(\cdot)$.\footnote{\PKL\ stands here for ``pi-kappa-lambda'' in the 
  sense of ``$\Pkl{}{}$''. } For 
a structure $\gmA=\pairof{A,\symb{K}^\gmA\ctenten}$, the weak second-order 
variables $X$, $Y$\ctenten\ run over elements of
$\Pkl{\symb{\scriptstyle K}^{\gmA}}{A}$.

We shall call a structure $\gmA$ with $\symb{K}$ in its signature 
as a unary predicate symbol \st\ $\cardof{\symb{K}^\gmA}$ is a regular 
uncountable cardinal, a {\it \PKL-structure}.

$\calL^{\PKL}_{stat}$ has 
the unique second-order quantifier ``$stat$'' and the internal interpretation $\models^{int}$ of 
formulas in this logic is defined similarly to $\calL^{\aleph_0}_{stat}$ 
with the crucial step in the inductive definition being
\begin{xitemize}
\xitem[PKL-0]   $\gmA\models^{\intnl}stat\,X\,\varphi(a_0\ctentenc U_0\ctentenc X)$\ \
  $\Leftrightarrow$\\
  \qquad$\setof{U\in\Pkl{\symb{\scriptstyle K}^\gmA}{A}\cap A}{\gmA\models^{\intnl}\varphi(a_0\ctentenc U_0\ctentenc U)}$
  is stationary in\\
  \qquad$\Pkl{\symb{\scriptstyle K}^\gmA}{A}$
\end{xitemize}
for an $\calL^{\PKL}_{stat}$-formula
$\varphi=\varphi(x_0\ctentenc X_0\ctentenc X)$ (for which the relation
$\models^{\intnl}$ has been defined), a \PKL-structure $\gmA=\pairof{A,\symb{K}^\gmA\ctenten}$ of a relevant signature,  
$a_0\ctenten\in A$ and $U_0\ctenten\in\Pkl{\symb{\scriptstyle K}^\gmA}{A}\cap A$.

For \PKL-structures $\gmA$, $\gmB$ of the same signature 
with $\gmB=\pairof{B,\symb{K}^\gmB\ctenten}$ and $\gmB\subseteq\gmA$, 
we define:
\begin{xitemize}
\xitem[PKL-1] 
  {\it $\gmB\prec^{\intnl}_{\calL^{\PKL}_{stat}}\gmA$\ \ $\Leftrightarrow$\\[2\jot]
  \qquad $\gmB\models^{\intnl}\varphi(b_0\ctentenc U_0\ctenten)$ if and only if 
  $\gmA\models^{\intnl}\varphi(b_0\ctentenc U_0\ctenten)$\\
  \qquad for all $\calL^{\PKL}_{stat}$-formulas $\varphi$ in the signature 
  of the structures with\\
  \qquad $\varphi=\varphi(x_0\ctentenc X_0\ctenten)$, $b_0\ctenten\in B$ and
  $U_0\ctenten\in\Pkl{\symb{\scriptstyle K}^\gmB}{B}\cap B$.}
\end{xitemize}

Finally, we define the internal strong downward L\"owenheim-Skolem theorem for this logic as 
follows:

Suppose that $\kappa$ is a regular cardinal $>\aleph_1$. 
\begin{xitemize}
\item[$\SDLS^{\intnl}_+(\calL^{\PKL}_{stat},\LT\kappa)$:]  {\it For any \PKL-structure
  $\gmA=\pairof{A,\symb{K}^\gmA\ctenten}$ of  
  countable signature with $\cardof{A}\geq\kappa$ and
  $\cardof{\symb{K}^\gmA}=\kappa$, there are {\it stationarily many}  
  $M\in[A]^{\LT\kappa}$ \st\  
  $\gmA\restr M$ is a \PKL-structure and
  $\gmA\restr M\prec^{\intnl}_{\calL^{\PKL}_{stat}}\gmA$.  }
\end{xitemize}

The following diagonal reflection characterizes
$\SDLS^{\intnl}_+(\calL^{\PKL}_{stat},\LT\kappa)$.
For regular cardinals $\kappa$, $\lambda$ with $\kappa\leq\lambda$, let
\begin{xitemize}
\item[$(\ast)^{\intnl+\PKL}_{\LT\kappa,\lambda}$: ] {\it For any countable expansion $\gmA$ 
  of the structure $\pairof{\calH(\lambda),\kappa,\in}$ and any family
  $\seqof{S_a}{a\in\calH(\lambda)}$ \st\ $S_a$ is a stationary subset of
  $\Pkl{\kappa}{\calH(\lambda)}$ for all $a\in\calH(\lambda)$, there are stationarily many 
  $M\in\Pkl{\kappa}{\calH(\lambda)}$ \st\ $\cardof{\kappa\cap M}$ is regular, $\gmA\restr M\prec\gmA$ and
  $S_a\cap \Pkl{\kappa\cap M}{M}\cap M$ is stationary in $\Pkl{\kappa\cap M}{M}$.}
\end{xitemize}

\begin{Proposition}\Label{P-PKL-0}{\rm (Proposition 4.1 in \cite{II})}
  For a regular cardinal $\kappa>\aleph_1$, \tfae:\smallskip

  \assert{a} $(\ast)^{\intnl+\PKL}_{\LT\kappa,\lambda}$ holds for all 
  regular $\lambda\geq\kappa$.\smallskip

  \assert{b} $\SDLS^{\intnl}_+(\calL^{\PKL}_{stat},\LT\kappa)$ holds.\qed
\end{Proposition}

For a regular cardinal $\kappa$ and a cardinal $\lambda\geq\kappa$, $\calS\subseteq\Pkl{}{}$ is 
said to be {\it $2$-stationary} if, for any stationary $\calT\subseteq\Pkl{}{}$, 
there is an $a\in\calS$ \st\ $\cardof{\kappa\cap a}$ is a regular uncountable 
cardinal and $\calT\cap\Pkl{\kappa\cap a}{a}$ is stationary in
$\Pkl{\kappa\cap a}{a}$. A regular cardinal $\kappa$ has the {\it$2$-stationarity 
  property} if $\Pkl{}{}$ is $2$-stationary (as a subset of itself) for all $\lambda\geq\kappa$. 

Since the property \assertof{a} in \Propof{P-PKL-0} is an extension of the 
$2$-stationarity of $\kappa$, we obtain: 

\begin{Lemma}
  \Label{P-PKL-1} For a regular uncountable $\kappa$,
  $\SDLS^{\intnl}_+(\calL^{\PKL}_{stat},\LT\kappa)$ implies that $\kappa$ is $2$-stationary.\qed
\end{Lemma}

This implies that a regular uncountable $\kappa$ with 
$\SDLS^{\intnl}_+(\calL^{\PKL}_{stat},\LT\kappa)$ must be a fairly large cardinal:

\begin{Lemma}\Label{P-PKL-2}{\rm (Lemma 4.3 in \cite{II})} Suppose that $\kappa$ is a regular uncountable cardinal.
  If $\kappa$ is $2$-stationary then $\kappa$ is a weakly Mahlo cardinal. \qed
\end{Lemma}

Actually the proof of \Lemmaof{P-PKL-2} (in \cite{II}) shows that $\kappa$ is weakly hyper 
Mahlo, weakly hyper hyper Mahlo, etc. 

\begin{Corollary}
  \Label{C-PKL-a}
  $\SDLS^{\intnl}_+(\calL^{\PKL}_{stat},\LT\continuum)$ implies that $\continuum$ is 
  weakly Mahlo, weakly hyper 
  Mahlo, weakly hyper hyper Mahlo, etc. \qed
\end{Corollary}

Using the characterization \Propof{P-PKL-0} of $\SDLS^{\intnl}_+(\calL^\PKL_{stat},\LT\kappa)$, 
a proof similar to that of \Thmof{T-refl-8} shows the following:
\begin{Theorem}
    \Label{P-lt-conti-3:} Suppose that $\kappa$ is a generically supercompact 
  cardinal by
  %% $\kappa$-cc
  $\LT\mu$-cc \pos\ for some $\mu<\kappa$. Then
  $\SDLS^{\intnl}_+(\calL^\PKL_{stat},\LT\kappa)$ holds. \qed
\end{Theorem}

Since ccc \pos\ are proper, we obtain the following 
by \Thmof{T-refl-10}, \Propof{P-internal-0} and \Thmof{P-lt-conti-3:}: 
\begin{Corollary}
  \Label{C-PKL-0}
  Suppose that $\kappa$ is a generically supercompact 
  cardinal by ccc \pos. Then $\SDLS^{\intnl}_+(\calL^{\aleph_0}_{stat},\LT\kappa)$ and
  $\SDLS^{\intnl}_+(\calL^\PKL_{stat},\LT\kappa)$ hold. \qed
\end{Corollary}
By \Lemmaof{L-PKL-0}, it follows that
\begin{Corollary}
  \Label{C-PKL-1} If \ZFC\ $+$ ``there is a supercompact cardinal'' is consistent 
  then so is \ZFC\ $+$ $\SDLS^{\intnl}_+(\calL^{\aleph_0}_{stat},\LT\continuum)$ and
  $\SDLS^{\intnl}_+(\calL^\PKL_{stat},\LT\continuum)$. Note that the continuum is 
  fairly large in the latter axiom system by \Corof{C-PKL-a}. \qed
\end{Corollary}
\section{Laver-generic large cardinals}
\Label{Laver}
The reflection properties we presented so far in connection with the size of the 
continuum can be summarized in three possible 
scenarios: \medskip

{\it
\wassert{A} \GRP\ \ \dotfill\ \ This implies \CH\ (\Lemmaof{L-refl-0}). \smallskip

\wassert{B} $\SDLS^-(\calL^{\aleph_0}_{stat},\LT\continuum)$\ \  \dotfill\ \ This implies
$\continuum=\aleph_2$ (\Thmof{T-refl-4},\,\assertof{2}).\smallskip

\wassert{C} $\SDLS^{\intnl}_+(\calL^{\PKL}_{stat},\LT\continuum)$\ \dotfill\ \mbox{This implies that 
$\continuum$ is fairly large  (\Corof{C-PKL-a}).}}\medskip

\ifextended{\tt
  \Thmof{T-refl-8} suggests that the existence of a generic large cardinal may be regarded 
  as a strong reflection principle. \assertof{A} above is actually such an instance 
  (see also \Corof{C-refl-1}). Thus, it seems to be a natural question whether the 
  scenarios \assertof{B} and \assertof{C} are also consequences of the existence of 
  certain generic large cardinals. 
}\fi %%\ifextended{tt

\ifextended{\tt Indeed, these}\else These\fi\ three possible scenarios can be treated in a uniform way from 
the point of view of \ifextended{\tt  the existence of}\fi\ the Laver-generic large cardinals defined below.

We shall call a class $\calP$ of \pos\  {\it iterable}  if
\begin{xitemize}
\xitem[laver-a-0] $\calP$ is closed \wrt\ forcing equivalence. That is, if 
  $\poP\in\calP$ and $\poP'$ is forcing equivalent to $\poP$ then $\poP'\in\calP$;
\xitem[laver-a-1] For any $\poP\in\calP$ and $\condp\in\poP$,
  $\poP\restr\condp\in\calP$;
\xitem[laver-a-2] If $\poP\in\calP$ and $\forces{\poP}{\utpoQ\varin\calP}$ then
  $\poP\ast\utpoQ\in\calP$. 
\end{xitemize}

Note that most of natural classes of \pos\ like $\sigma$-closed \pos, ccc \pos, proper \pos, stationary 
preserving \pos\ etc. are iterable.

For a cardinal $\kappa$ and an iterable class $\calP$ of \pos, we call $\kappa$ a 
{\it Laver-generically supercompact for $\calP$} if, for any $\lambda\geq\kappa$ and any
$\poP\in\calP$, there is a \po\ $\poQ\in\calP$ 
with $\poP\circleq\poQ$ \st, for any $(\uniV,\poQ)$-generic filter $\genH$,  there 
are $M$, $j\subseteq\uniV[\genH]$ \st\ 
\begin{xitemize}
\xitem[laver-a-3] $M$ is an inner model of $\uniV[\genH]$, 
\xitem[laver-a] $\elembed{j}{\uniV}{M}$, 
\xitem[laver-0] $\crit(j)=\kappa$, $j(\kappa)%% =\cardof{\poQ}
  >\lambda$,
\xitem[laver-0-0] $\poP$, $\genH\in M$ and 
\xitem[laver-1] $j\imageof\lambda\in M$.\footnote{This definition of 
  Laver-generic supercompactness for $\calP$ is different  
  from the one given in \cite{II}. However, it is easy to show that the present 
  definition is equivalent to the one in \cite{II} for an iterable $\calP$. Note 
  that, strictly speaking, this equivalence is used at the end of the proof of 
  \Thmof{T-laver-4} below.}
\end{xitemize}

$\kappa$ is {\it Laver-generically superhuge} ({\it Laver-generically super almost-huge} 
resp.) {\it for $\calP$} if $\kappa$ satisfies 
the definition of Laver-generic supercompactness for $\calP$ with  
\xitemof{laver-1} replaced by 
\begin{xitemize}
\xitemd[laver-1]{'}
  $j\imageof j(\kappa)\in M$  ($j\imageof\mu\in M$ for all $\mu<j(\kappa)$ resp.).
\end{xitemize}

$\kappa$ is {\it tightly Laver-generically supercompact} ({\it tightly 
  Laver-generically superhuge, tightly Laver-generically super almost-huge}, 
resp.) if the definition of {Laver-generically supercompact} ({Laver-generically 
  superhuge, Laver-generically super almost-huge}, resp.) holds with 
\xitemof{laver-0} replaced by 
\begin{xitemize}
\xitemd[laver-0]{'} 
  $\crit(j)=\kappa$, $j(\kappa)=\cardof{\poQ}>\lambda$.
\end{xitemize}

All consistency proofs of the existence of Laver-generic large cardinals we know actually 
show the existence of tightly Laver-generic large cardinals (see 
the proof of Theorem 5.2 in \cite{II}). 

The consistency of the existence of a Laver-generic large cardinal can be proved from the 
assumption of the existence of the corresponding genuine large cardinals except the case of 
the Laver-generic large cardinals by proper \pos. This case will be further 
discussed in \cite{III}. 

\begin{Theorem}\Label{T-laver-0}{\rm (\cite{II})}
    \wassertof{1} Suppose that \ZFC\ $+$ ``there exists a supercompact cardinal'' is 
  consistent. Then \ZFC\ $+$ ``there exists a Laver-generically supercompact cardinal 
  for 
  $\sigma$-closed \pos'' is consistent as well.\smallskip
  
  \wassert{2} Suppose that \ZFC\ $+$ ``there exists a superhuge cardinal'' is 
  consistent. Then \ZFC\ $+$ ``there exists a Laver-generically super almost-huge 
  cardinal for proper \pos'' is consistent as well.\smallskip

  \wassert{3} Suppose that \ZFC\ $+$ ``there exists a supercompact cardinal'' is 
  consistent. Then \ZFC\ $+$ ``there exists a Laver-generically supercompact 
  cardinal $\kappa$ for c.c.c.\ \pos'' with $\kappa=2^{\aleph_0}$ is consistent as well.
\end{Theorem}
{\bf Sketch of the proof.} Let us consider \assertof{2}. The other assertions are 
similarly and easier to prove. 

Starting from a model of \ZFC\ with a superhuge cardinal $\kappa$, we can 
obtain models of respective assertions by 
iterating in countable support with proper \pos\ $\kappa$ times along a Laver function
for super almost-hugeness which exists by a result in Corazza \cite{corazza}.

In the resulting model, we obtain Laver-generically super almost-hugeness in 
terms of proper \po\ $\poQ$ in each respective inner model $M[\genG]$ of
$\uniV[\genG]$. The closedness of $M$ in $\uniV$ in terms of super almost-hugeness 
implies that $\poQ$ is also proper in $\uniV[\genG]$ (this is the place where we 
need the super almost-hugeness: for \assertof{1} and \assertof{3} we do not need 
this much closedness of $M$).

This shows that $\kappa$ 
is Laver-generically super almost-huge of proper \pos. 
\qedofThm\qedskip

In contrast to simple generic supercompactness, a Laver-generically 
supercompact cardinal for a natural class $\calP$ of \pos\ is determined 
uniquely if it exists:

\begin{Proposition}\Label{P-laver-0}{\rm (\cite{II})}
  \wassertof{0} If $\kappa$ is generically measurable for some \po\ $\poP$, then 
  $\kappa$ is regular.\smallskip

  \wassert{1}   Suppose that $\kappa$ is generically measurable by a $\omega_1$ 
  preserving $\poP$. Then $\kappa>\omega_1$.\smallskip

  \wassert{2}  Suppose that $\kappa$ is Laver-generically supercompact for
  a class $\calP$ of \pos\ with $\Col(\omega_1,\ssetof{\omega_2})\in\calP$. 
  Then $\kappa\leq\omega_2$. If all elements of $\calP$ are $\omega_1$-preserving, 
  then we have $\kappa=\omega_2$. \smallskip

  \wassert{3} Suppose that $\calP$ is a class of \pos\ containing 
  a \po\ $\poP$ \st\ any $(\uniV,\poP)$-generic filter $\genG$ codes a new real. If $\kappa$ is a 
  Laver-generically supercompact for $\calP$, 
  then $\kappa\leq\continuum$.\smallskip

  \wassert{4} Suppose that $\calP$ is a class of \pos\ \st\ elements of 
  $\calP$ do not add any reals. If $\kappa$ is Laver-generically supercompact by $\calP$, then 
  we have $\continuum<\kappa$.%% \smallskip
  %%
  %% \wassert{5} Suppose that $\kappa$ is generically supercompact by $\calP$ \st\ there is at least one 
  %% $\poP\in\calP$ which adds a new real. Then $\kappa\leq 2^{\aleph_0}$ holds.
  \ifextended\else\qed\fi
\end{Proposition}
\ifextended{\tt \prf \assertof{0}: Suppose not and let $\mapping{f}{\mu}{\kappa}$ 
 for some  $\mu<\kappa$ be cofinal. Let $\poP$, $\genG$, $j$, $M$ be as in the 
 definition of generic measurability of $\kappa$. Since $j(f)=f$ and
 $\uniV\models\kappa\equiv\sup(f)$. By elementarity, it follows that
 $M\modelof{j(\kappa)\equiv\sup\underbrace{j(f)}_{=f}}$. Thus $j(\kappa)=\kappa$. 
 This is a contradiction.\smallskip

 \assertof{1}: Suppose not. Then $\kappa=\omega_1$. Let $\genG$ 
 be $(\uniV,\poP)$-generic and $\elembed{j}{\uniV}{M\subseteq\uniV[\genG]}$ be 
 \st\ $\crit(j)=\omega_1$. By elementarity, $M\modelof{j(\kappa)\equiv\omega_1}$. 
 Thus $M\modelof{\kappa\xmbox{ is countable}}$. It follows that
 $\uniV[\genG]\modelof{\kappa\xmbox{ is countable}}$. This is a contradiction to 
 the $\omega_1$-preserving of $\poP$. \smallskip

 \assertof{2}: Suppose $\kappa>\omega_2$. Let $\poQ\in\calP$ be \st\ $\Col(\omega_1,\ssetof{\omega_2})\circleq\poQ$
 be with $(\uniV,\poQ)$-generic filter $\genH$ and $M$, $j\subseteq\uniV[\genH]$ 
 \st\ $\elembed{j}{\uniV}{M}$, $\crit(j)=\kappa$ and
 $\Col(\omega_1,\ssetof{\omega_2})$, $\calH\in M$. Then $j(\omega_2)=\omega_2$ 
 and thus $M\modelof{(\omega_2)^\uniV\mbox{ is }\omega_2}$. But
 $\Col(\omega_1,\ssetof{\omega_2})$ part of $\genH$ in $M$ collapses
 $(\omega_2)^\uniV$ to an ordinal $\LE\omega_1$. This is a contradiction which shows 
 that $\kappa\leq\omega_2$. 

 If all the elements of $\poP$ is $\omega_1$-preserving we also obtain
 $\kappa\geq\omega_2$ by \assertof{1} and hence $\kappa=\omega_2$.
 \smallskip

 \assertof{3}: Suppose that $\mu<\kappa$ and $\seqof{a_\xi}{\xi<\mu}$ is a 
 sequence of reals. It is enough to show that this sequence does not enumerate 
 all reals. Let $\poQ\in\calP$ be \st\ $\poP\circleq\poQ$ 
 with $(\uniV,\poQ)$-generic filter $\genH$ and $M$, $j\subseteq\uniV[\genH]$ 
 \st\ $\elembed{j}{\uniV}{M}$, $\crit(j)=\kappa$ and
 $\Col(\omega_1,\ssetof{\omega_2})$, $\genH\in M$. By $\mu<\crit(j)$ we have
 $j(\seqof{a_\xi}{\xi<\mu})=\seqof{a_\xi}{\xi<\mu}$. Since $M$ contains a new 
 real added by $\poP$,
 $M\modelof{\seqof{a_\xi}{\xi<\mu}\xmbox{ is not an enumeration of all reals}}$. 
 By elementarity it follows that 
 $\uniV\modelof{\seqof{a_\xi}{\xi<\mu}\xmbox{ is not an enumeration of all reals}}$. 
 \smallskip

 \assertof{4}: Suppose $\kappa\leq\continuum$. Let $\lambda>\continuum$ and let
 $\poP\in\calP$ be with $(\uniV,\poP)$-generic $\genG$ and $M$, 
 $j\subseteq\uniV[\genG]$ \st\ $\elembed{j}{\uniV}{M}$, $\crit(j)=\kappa$ and
 $j(\kappa)>\lambda$. Since $\poP$ does not add any new real, it follows that 
 \begin{xitemize}
 \xitema[] $M\modelof{\underbrace{\continuum}_{=(\continuum)^V}\geq j(\kappa)>\lambda>(\continuum)^V}$.
 \end{xitemize}
 A contradiction. 
\qedofProp\qedskip
}\fi %\ifextended{\tt 

For a class $\calP$ of \pos\ and cardinals $\mu$, $\kappa$, 

\begin{xitemize}
\item[$\MA^{+\mu}(\calP,\LT\kappa)$:]
  \it For any $\poP\in\calP$, any family 
  $\calD$ of dense subsets of\/ $\poP$ with $\cardof{\calD}<\kappa$ and any family 
  $\calS$ of\/ $\poP$-names \st\ $\cardof{\calS}\leq\mu$ and
  $\forces{\poP}{\utilde{S}\xmbox{ is a stationary subset of }\omega_1}$ for all
  $\utilde{S}\in\calS$, there is a $\calD$-generic filter $\genG$ over $\poP$ \st\
  $\utilde{S}[\genG]$ is a stationary subset of $\omega_1$ for all
  $\utilde{S}\in\calS$. 
\end{xitemize}

The following strengthening of the Laver-genericity is needed to obtain  
``$++$'' versions of forcing axioms. 

For a cardinal $\kappa$ and an iterable class $\calP$ of \pos, we call $\kappa$ a 
{\it strongly Laver-generically supercompact for $\calP$} if, 
for any $\lambda\geq\kappa$ and any 
$\poP\in\calP$, there is a $\poP$-name of a \po\ $\utpoQ$ with $\forces{\poP}{\utpoQ\varin\calP}$
\vspace{\smallskipamount}\st, for any $(\uniV,\poP\ast\utpoQ)$-generic filter $\genH$,  there 
are $M$, $j\subseteq\uniV[\genH]$ with 
\begin{xitemize}
\xitem[laver+-a-3] $M$ is an inner model of $\uniV[\genH]$, 
\xitem[laver+-a] $\elembed{j}{\uniV}{M}$, 
\xitem[laver+-0] $\crit(j)=\kappa$, $j(\kappa)%% =\cardof{\poQ}
  >\lambda$,
\xitem[laver+-0-0] $\poP$, $\genH\in M$ and 
\xitem[laver+-1] $(\fnsp{\lambda}{M})^{\uniV[\genH]}\subseteq M$.
\end{xitemize}
The notions of {\it strongly Laver-generically superhuge} and {\it strongly Laver-generically super almost-huge} 
are defined correspondingly.
For example, $\kappa$ is strongly Laver-generically super almost-huge if the definition 
of strongly Laver-generically supercompact cardinal holds with \xitemof{laver+-1} 
replaced with
\begin{xitemize}
\xitem[laver+-2] $(\fnsp{\mu}{M})^{\uniV[\genH]}\subseteq M$ holds for all $\mu<j(\kappa)$. 
\end{xitemize}

Note that, if $\calP$ is the class of ccc \pos, then the ``strongly'' version of the 
Laver-generically large cardinal is equivalent to the original version of the 
corresponding Laver-generic largeness. 
Note also that the construction in the proof of \Thmof{T-laver-0} actually 
provides models of strongly Laver-genericity. 

\begin{Theorem}\Label{T-laver-1}{\rm(Theorem 5.7 in \cite{II})}
  For a class $\calP$ of proper 
  \pos, if $\kappa>\aleph_1$ is a strongly Laver-generically supercompact for $\calP$, 
  then $\MA^{+\mu}(\calP,\LT\kappa)$ holds for all $\mu<\kappa$.\qed
\end{Theorem}

\begin{Lemma}
  \Label{L-laver-0}
  Suppose that $\kappa$ is generically supercompact by a class $\calP$ of 
  \pos\ \st\ all $\poP\in\calP$ has the $\LT\mu$-cc for some $\mu<\kappa$. Then for
  any regular $\lambda\geq\kappa$, $\Pkl{}{}$ has a $\mu$-saturated normal fine 
  filter $F$ over $\Pkl{}{}$.
\end{Lemma}
\prf For a regular $\lambda\geq\kappa$, let $\poP$ be \st\ there are 
$(\uniV,\poP)$-generic $\genG$, and $M$, $j\subseteq\uniV[\genG]$ \st\ 
$\elembed{j}{\uniV}{M}$, $\crit(j)=\kappa$, $j(\kappa)>\lambda$ and
$j\imageof\lambda\in M$. \vspace{\jot}

In $\uniV$, 
let $F=\setof{A\subseteq\Pkl{}{}}{\forces{\poP}{j\imageof\lambda\varin j(\check{A})}}$. 
By the $\LT\mu$-cc of $\poP$, this $F$ is as desired.\qedofLemma
\qedskip

Combining \Propof{P-laver-0}, \Thmof{T-laver-1} and \Lemmaof{L-laver-0}, we 
obtain:

\begin{Theorem}\Label{T-laver-2} Suppose that $\kappa$ is 
  strongly Laver-generically supercompact cardinal for an iterable class $\calP$ of 
  \pos.\medskip

  \wassert{A$'$} If all elements of $\calP$ are $\omega_1$-preserving and do not add any 
reals, and\\ $\Col(\omega_1,\ssetof{\omega_2})\in\calP$, 
then $\kappa=\aleph_2$ and \CH\ holds.\smallskip

\wassert{B$'$} If all elements of $\calP$ are $\omega_1$-preserving and $\calP$ contains all 
proper \pos\ then $\PFA^{+\omega_1}$ holds and $\kappa=\continuum=\aleph_2$.
\smallskip

\wassert{C$'$} If all elements of $\calP$ are $\LT\mu$-cc for some fixed $\mu<\kappa$ and 
$\poP$ contains a \po\ which adds a new real then $\kappa$ is fairly large (in the 
sense of \Lemmaof{L-laver-0}) and $\kappa\leq\continuum$. 
\end{Theorem}

\prf \assertof{A$'$}\,: By \Propof{P-laver-0},
\assertof{2} and \assertof{4}.\smallskip

\assertof{B$'$}\,: $\kappa=\aleph_2$ by \Propof{P-laver-0}, \assertof{2}. 
$\PFA^{+\omega_1}$ holds by \Thmof{T-laver-1}. $\PFA$ 
implies $\continuum=\aleph_2$. \smallskip

\assertof{C$'$}\,: $\kappa\leq\continuum$ by \Propof{P-laver-0}, \assertof{3}. 
$\kappa$ is large by \Lemmaof{L-laver-0}.\qedofThm 
\qedskip

The three cases in \Thmof{T-laver-2} can be further modified to fit to the reflection 
principles discussed in earlier sections.

\begin{Theorem}\Label{T-laver-3} \wassertof{A$''$} Suppose that $\kappa$ is Laver-generically 
  supercompact for $\sigma$-closed \pos. Then $2^{\aleph_0}=\aleph_1$,
  $\kappa=\aleph_2$, $\MA^{+\omega_1}(\sigma\mbox{-closed})$ and 
  \GRP\ holds. It follows that \RC\ and
  $\SDLS(\calL^{\aleph_0}_{stat},\LT\aleph_2)$ hold. \smallskip

  \wassert{B$''$} Suppose that elements of $\calP$ are $\omega_1$-preserving and $\calP$ contains all 
  proper \pos. If $\kappa$ is strongly Laver-generically supercompact 
  for $\calP$, then 
  $2^{\aleph_0}=\kappa=\aleph_2$, $\PFA^{+\omega_1}$ and hence also 
  $\SDLS^-(\calL^{\aleph_0}_{stat},\LT\continuum)$ holds.\smallskip

  \wassert{C$''$} Suppose that $\kappa$ is Laver-generically supercompact for ccc 
  \pos. Then $\kappa\leq2^{\aleph_0}$ and $\Pkl{}{}$ for any regular $\lambda\geq\kappa$ carries an
  $\aleph_1$-saturated normal ideal. In particular, $\kappa$ is $\kappa$-weakly Mahlo.
  $\MA^{+\mu}(ccc,\LT\kappa)$ for all $\mu<\kappa$,  
  $\SDLS^{\intnl}(\calL^{\aleph_0}_{stat},\LT\kappa)$ and 
  $\SDLS^{\intnl}_+(\calL^{\PKL}_{stat},\LT\kappa)$ also hold.
\end{Theorem}
\prf \assertof{A$''$}\ : $2^{\aleph_0}=\aleph_1$ and $\kappa=\aleph_2$ follows from 
\Thmof{T-laver-2},\,\assertof{A}. $\MA^{+\omega_1}(\sigma\mbox{-closed})$ holds 
by \Thmof{T-laver-1}. \GRP\ holds by \Corof{C-refl-1}. \RC\ and
  $\SDLS(\calL^{\aleph_0}_{stat},\LT\aleph_2)$ follow from \GRP\ by 
\Thmof{T-refl-7} and \Corof{C-refl-0},\,\assertof{1}. 
\smallskip

\assertof{B$''$}\ : This is just as \assertof{B} in \Thmof{T-laver-2}. 
$\SDLS^-(\calL^{\aleph_0}_{stat},\LT\continuum)$ holds by 
\Corof{C-refl-a},\,\assertof{1}. 
\smallskip

\assertof{C$''$}\ : The first half of the assertion follows from 
\Thmof{T-laver-2},\,\assertof{C}. \\
$\MA^{+\mu}(ccc,\LT\kappa)$ for 
all $\mu<\kappa$ holds by \Thmof{T-laver-1}. 
$\SDLS^{\intnl}(\calL^{\aleph_0}_{stat},\LT\kappa)$ and 
$\SDLS^{\intnl}_+(\calL^{\PKL}_{stat},\LT\kappa)$ hold by \Corof{C-PKL-0}.
\qedofThm
\qedskip

At the moment we do not know whether the assumption in \assertof{C} in 
\Thmof{T-laver-3} implies $\kappa=\continuum$. The following partial 
answer is obtained in \cite{II}:
\begin{Theorem}\Label{T-laver-3-0}{\rm(Theorem 5.3 in \cite{II})}
  If $\kappa$ is tightly Laver-generically superhuge for ccc \pos, then
  $\kappa=\continuum$. \qed
\end{Theorem}

The following \Thmof{T-laver-4} is a Laver-generic version of \Thmof{T-refl-9}. The proof 
of the theorem is a typical application of the master condition argument (see 
e.g.\, Cummings \cite{cummings}). 

\begin{Theorem}
  \Label{T-laver-4} Suppose that $\kappa$ is strongly Laver-generically almost super-huge for 
  an iterable $\calP$ which provably contains all \pos\ of the form $\Col(\mu,\lambda)$ 
  for all regular $\mu$, $\lambda$ with $\kappa\leq\mu<\lambda$. For a regular
  $\lambda_0$, let $\poP_0=\Col(\kappa,\lambda_0)$ and let $\genG_0$ be a
  $(\uniV,\poP_0)$-generic filter.

  Then we have
  \begin{xitemize}
  \xitem[laver-2] 
    $\uniV[\genG_0]\modelof{{}
    \begin{array}[t]{@{}l}
      \kappa\mbox{ is strongly Laver-generically super almost-huge for }\calP}.
    \end{array}
    $
  \end{xitemize}
\end{Theorem}
\prf %%In the following proof, we are using Lemma 2.5 in \cite{II} freely without mention. 
Suppose that $\uniV[\genG_0]\modelof{\poP\in\calP}$ and let $\utpoP$ be a 
$\poP_0$-name of $\poP$ \st\ $\forces{\poP_0}{\utpoP\varin\calP}$. Let $\lambda\geq\kappa$ be regular. \Wolog, we may 
assume that $\lambda\geq\cardof{\poP}\geq\lambda_0$. 
$\poP_0\ast\utpoP\in\calP$ by iterability of $\calP$. 
Since $\kappa$ is strongly Laver-generically super almost-huge for $\calP$, 
there are 
$\poP_0\ast\utpoP$-name $\utpoQ$ with 
$\forces{\poP_0\ast\utpoP}{\utpoQ\varin\calP}$ and $(\uniV,\poP_0\ast\utpoP\ast\utpoQ)$-generic 
filter $\genH$ \st\ $\genG_0\subseteq\genH$ and \st\ there are $M$,
$j\subseteq\uniV[\genH]$ with \xitemof{laver+-a-3}, \xitemof{laver+-a}, 
\xitemof{laver+-0}, \xitemof{laver+-2} and
\begin{xitemize}
\xitem[laver-3] $\poP_0\ast\utpoP$, $\genH\in M$. 
\end{xitemize}
We have $j\imageof\genG_0\in M$. Let $\poP_1=j(\poP_0)$ by \xitemof{laver+-2}.

By elementarity, we have
$M\modelof{\poP_1=\Col(j(\kappa),j(\lambda_0))}$. Note that we also have 
$\uniV[\genH]\modelof{\poP_1=\Col(j(\kappa),j(\lambda_0))}$ by 
\xitemof{laver+-2}. 
Since
$M\modelof{j\imageof\genG_0\xmbox{ has the fip}}$ by elementarity and 
$M\modelof{\cardof{j\imageof\genG_0}\leq\cardof{\poP}\leq\lambda<j(\kappa)}$, 
there is $\condq^*\in \poP_1$ in $M$ \st\ $\condq^*\leq_{\poP_1}j(\condp)$ for all $\condp\in\genG_0$. 
Let $\genG_1$ be $(\uniV[\genH],\poP_1)$-generic filter with $\condq^*\in\genG_1$. 
In $\uniV[\genH\ast\genG_1]$, let 
\begin{xitemize}
\xitem[laver-4] 
  $\elembed{\tilde{j}}{\uniV[\genG_0]}{M[\genG_1]}$;
  $\utilde{a}[\genG_0]\mapsto j(\utilde{a})[\genG_1]$. 
\end{xitemize}
Since $\poP\in\calP^{\uniV[\genG_0]}$ and $\lambda$ were arbitrary, the 
elementary embedding $\tilde{j}$ above witnesses the Laver-generic 
super almost-hugeness of $\kappa$ for $\calP$ in $\uniV[\genG_0]$. 
\qedofThm

\begin{Corollary}
  \Label{C-laver-0} Suppose that 
$\calP$ is an iterable class of \pos\ which provably contains all \pos\ of the form $\Col(\mu,\lambda)$ 
  for all regular uncountable $\mu$, $\lambda$ with $\kappa\leq\mu<\lambda$. If the theory \ZFC\ $+$ 
  ``there is a strongly Laver-generically super almost-huge $\kappa$ for $\calP$'' $+$ ``there is a 
  supercompact $\kappa_1>\kappa$'' is consistent, then so is the theory \ZFC\ $+$ ``there is 
  a strongly Laver-generically super almost-huge $\kappa$ for $\calP$'' $+$
  ``$\kappa^+$ is generically supercompact by $\LT\kappa$-closed \pos''. In 
  particular, $\GRP^{\LT\kappa}(\LT\kappa^+)$ follows from this theory.  
\end{Corollary}
\prf Suppose that $\kappa$ is strongly Laver-generically supercompact for $\calP$ and 
$\kappa_1$ is a supercompact cardinal. Let $\poP_0=\Col(\kappa,\kappa_1)$ and 
let $\genG_0$ be a $(\uniV,\poP_0)$-generic filter. By \Thmof{T-laver-4},
$\uniV[\genG_0]\modelof{\kappa\xmbox{ is strongly Laver-generically supercompact for }\calP}$. 
$\uniV[\genG_0]\models\kappa_1=\kappa^+$ and
$\uniV[\genG_0]\modelof{\kappa^+\xmbox{ is generically supercompact by }\LT\kappa\xmbox{-closed 
\pos}}$ by \Lemmaof{L-refl-1}. By 
\Thmof{T-refl-8}, $\GRP^{\LT\kappa}(\LT\kappa^+)$ follows.\qedofCor
\qedskip

\Thmof{T-laver-4} and \Corof{C-laver-0} have many variants with similar proofs. 
For example:
\begin{Theorem}
  \Label{T-laver-5} Suppose that $\kappa$ is strongly Laver-generically supercompact for 
  an iterable $\calP$ which provably contains all $\sigma$-closed \pos. For a regular
  $\lambda_0$, let $\poP_0=\Col(\kappa,\lambda_0)$ and let $\genG_0$ be a
  $(\uniV,\poP_0)$-generic filter.

  Then we have
  \begin{xitemize}
  \xitem[laver-5] 
    $\uniV[\genG_0]\modelof{{}
    \begin{array}[t]{@{}l}
      \kappa\mbox{ is strongly Laver-generically supercompact for }\calP}. 
    \end{array}
    $ \qed
  \end{xitemize}
\end{Theorem}

\begin{Corollary}
  \Label{C-laver-1} Suppose that 
  $\calP$ is an iterable class of \pos\ which provably contains all $\sigma$-closed \pos. If the theory \ZFC\ $+$ 
  ``there is a strongly Laver-generically supercompact $\kappa$ for $\calP$'' $+$ ``there is a 
  supercompact $\kappa_1>\kappa$'' is consistent, then so is the theory \ZFC\ $+$ ``there is 
  a strongly Laver-generically supercompact $\kappa$ for $\calP$'' $+$
  ``$\kappa^+$ is generically supercompact by $\LT\kappa$-closed \pos''. In 
  particular, $\GRP^{\LT\kappa}(\LT\kappa^+)$ follows from this theory.  \qed
\end{Corollary}

Note that, by \Propof{P-laver-0},\,\assertof{2}, we have $\kappa\leq\aleph_2$ in 
\Thmof{T-laver-5} and \Corof{C-laver-1} above. 

In \cite{III}, we show that the combination of the principles $\SDLS^{int}(\calL^{\aleph_0}_{stat},\LT\continuum)$ and
$\GRP^{\LT\continuum}(\LE\continuum)$ is also consistent under large continuum 
assuming the consistency of two supercompact cardinals. 

\section{Some open problems and Takeuti's account on G\"odel's contribution to 
  the continuum problem}
\Label{epilogue}
Let us mention some open problems. Some of them will be addressed in 
\cite{III}. 

The following problem is already mentioned in the previous section:
\begin{problem}
  If $\kappa$ is Laver-generically supercompact for ccc \pos, does this imply 
  $\kappa=\continuum$? 
\end{problem}

Forcing axioms have some characterizations which may be interpreted as 
suggestions of the correctness of the axioms. See 
e.g. Bagaria \cite{bagaria0}, Fuchino \cite{potential}. 

\begin{problem}
  Is there any nice characterizations of ``$+\mu$'' versions of forcing axioms? 
\end{problem}
Any meaningful answer to this problem would enhance the relevance of the 
trichotomy in \Thmof{T-laver-3}. 

The trichotomy \assertof{A}, \assertof{B}, \assertof{C} of reflection principles 
mentioned at the beginning of \sectionof{Laver} 
has an alternative trichotomy \assertof{A}, \assertof{B${}^*$}\,, 
\assertof{C${}^*$}\,, where 
\medskip

{\it
\wassert{B${}^*$} $\SDLS^-(\calL^{\aleph_0}_{stat},\LT\continuum)$ and $\RC$ hold.
\smallskip

\wassert{C${}^*$} $\SDLS^{\intnl}_+(\calL^{\PKL}_{stat},\LT\continuum)$ and
$\RC(\LT\continuum)$ hold.
}\medskip

Note that \RC\ in \assertof{B${}^*$} is equivalent to $\RC(\LT\continuum)$ since
$\SDLS^-(\calL^{\aleph_0}_{stat},\LT\continuum)$ implies $\continuum=\aleph_2$.

\assertof{B${}^*$} and \assertof{C${}^*$} are not compatible with 
\assertof{B$''$} and \assertof{C$''$} respectively, since $\MA(\kappa)$ --- i.e. \MA\ (for ccc \pos) 
for $\LE\kappa$ many dense sets --- implies the negation of $\RC(\LE\kappa^+)$ 
(see section 5 in Fuchino \cite{fuchino3}). 

\assertof{B${}^*$} and \assertof{C${}^*$} can be realized by starting from a 
supercompact cardinal and then forcing with a Mitchell type mixed support 
iteration (for \assertof{B${}^*$}\,, this is mentioned in \cite{todorcevic}, for 
\assertof{C${}^*$}\,, see \cite{III}). Actually, \assertof{C${}^*$} can be realized 
much easier by simply adding supercompact many Cohen reals, but we do need mixed support 
iteration to obtain a model of\medskip

{\it
\wassert{C${}^\dagger$} $\SDLS^{\intnl}_+(\calL^{\PKL}_{stat},\LT\continuum)$, 
$\RC(\LT\continuum)$ and $\GRP^{\LT\continuum}(\LE\continuum)$
hold.
}\medskip
\medskip

Models obtained in this way seem to be much more artificial than the models for 
\assertof{B$''$} and \assertof{C$''$} as in the proof of \Thmof{T-laver-0}. Even so, we have the feeling that 
we do not understand yet much about the models constructed by the mixed support 
iteration and its variations. Thus
\begin{problem}\Label{mixed-support}
  What is possible with variations of mixed support iteration? 
\end{problem}

The result of K\"onig in \cite{koenig} cited here as \Corof{C-refl-1} can be also regarded as a 
characterization of $\aleph_2$ being generically supercompact by $\sigma$-closed 
\pos. Thus we may further ask:
\begin{problem}
  Is there any characterizations of Laver-generically large cardinals in terms of 
  some strong reflection principles? 
\end{problem}

The following essential problem might be much harder than the other problems: 
\begin{problem}
  Are Laver-generically large cardinals equiconsistent with corresponding genuine 
  large cardinals?
\end{problem}

The consistency of the reflection of non-metrizability of a first countable topological spaces 
down to $\LT\aleph_2$ is an open problem known as Hamburger's problem.
The consistency of the reflection of the property of partial orderings that they are not represented 
as countable union of chains down to $\LT\aleph_2$ is also an open problem known 
as Galvin's conjecture. In case of Hamburger's problem, it is known that the 
reflection of non-metrizability of a first countable topological spaces down to
$\LT\continuum$ is realized in the model obtained by adding supercompact many 
Cohen reals (Dow, Tall and Weiss \cite{dow-tall-weiss}). The reflection cardinal 
of uncountable chromatic number of graphs is known to be $\geq\beth_\omega$ 
(Erd\H{o}s and Hajnal, see \cite{fuchino-sakai} for a detailed proof in \ZFC). 
There are many open problems in connection with reflection of these and some other 
mathematical properties. Some of them seem to be extremely difficult. 
Let us mention here merely one problem which may have some connection to \Problemof{mixed-support}:
\begin{problem}
  Is the reflection of non-metrizability of first countable topological spaces 
  down to $\LT\continuum$ consistent with $\RC(\LT\continuum)$? 
\end{problem}

The first author of this article belongs to the generation of Japanese logicians 
who were strongly inspired by the writings of late Professor Gaishi Takeuti who 
published many expository articles and books in Japanese from 1960's to the end 
of 1990's. Although set theory was not his main field, Professor Takeuti wrote many 
expositions and told his views on the subject.  G\"odel's program was 
one of the issues he discussed repeatedly there.

G\"odel, who usually refused to publish papers which 
he thought was not yet perfect,  tried once to publish a quite unfinished note in 
1970 with the title ``Some considerations leading to the probable conclusion that the 
true power of the continuum is $\aleph_2$.'' It is said that, being seriously ill, he did so under 
the fear that he would soon die. After this crisis, G\"odel withdrew the note 
finding out some inaccuracy in it but he continued the study on the problem.
Oskar Morgenstern noticed in 
his diary on 20. September, 1975 that G\"odel told him in a telephone call that 
he was finally convinced that (the  
newest version of) his axiom implies that the continuum is ``different from
$\aleph_1$'' and that he will write it up (\cite{morgenstern}). Dawson 
\cite{dawson} contains some accounts about this development. 

\iffalse
Sa. 20. Sept.

Noch immer das grässliche Wetter: trüb, Tiefdruck, keine Luft, feucht – Das trägt nicht zum Wohlbefinden bei.
Mit Gerry tel. Es wird alles klappen. Er ist sehr fleissig gewesen.
Endlich Brief von Schwödianer: Sie haben einen Sohn, er war in Stockholm & 
Toronto. In T. gab er "unser" paper, wird es schicken (ha!) & tel. nicht von T., 
was so einfach gewesen wäre! Er ist das 2. Sem. in NYU. Im Jan. eine 
Kuratoriumssitz. in Wien. – werde ich das noch schaffen?? Noch später: Eben rief 
Jim Varney an: IVP war OK – aber unterhalb der Oper.stelle von vor 1 Jahr (heute, 
ausgerechnet!) eine Stelle: Krebs. Ich solle 5-10 Xray treetments bekommen; dann 
würden die Schmerzen nachlassen & man werde den Tumor stoppen (?). So zu Dr 
Fine & dann 5-10 mal je eine Dose. 
Später:
Wir hatten Gödel Blumen geschickt mit einer Gratulation zur Pres. Medal. Er rief 
bald an: "bubbling over". Wie gross & schwer sie sei. viell. Gold?" … und noch 
eine kleine für das Knopfloch" … Ich freute mich über alles, wie auch Dorothy, 
die er zuerst erreichte. 
Dann: jetzt weiss er daß sein Beweis über die Mächtigkeit des Cont.n. richtig 
ist. Es ist nicht Aleph 2, sondern: "verschieden von Aleph 1". Jetzt wird er das 
wirklich aufschreiben. Wir wollen uns bald treffen, viell. in 8 Tagen, ein 
Empfang im Institut. 
\fi

Takeuti wrote about the details of what he 
understood from the 1970 note on pp.99--124 in his book \cite{takeuti1} published 
in 1972. According to 
Takeuti \cite{takeuti3}, he was then invited by G\"odel in 1975 or 
1976\ifextended\footnote{\tt According to \cite{takeuti4}, it was 1976 since 
  Takeuti writes there that it was about a month before G\"odel's retirement.}\fi\ 
shortly 
before G\"odel's retirement from IAS and 
discussed with him about the results on the Continuum Problem. Takeuti 
\cite{takeuti2} in 1978 must be closely related to this discussion.

Twenty years later, in May 1998, the first author of the present article obtained a letter from Professor 
Takeuti with a copy of his handwritten manuscript in Japanese, which contained a further development of the 
material in \cite{takeuti2} among other things. 
Unfortunately, the first author could not give any  
reasonable comments to the manuscript at that time. 
The part of the manuscript on ``G\"odel's Continuum Hypothesis'' was then published 
in the new edition of \cite{takeuti3} as an appendix in September 1998. 

\iffalse{\tt
  There is a popular exposition by J.\ Kennedy \cite{kennedy} on this 
  development\footnote{\tt \cite{kennedy} is a nicely written exposition for a 
    general audience. However, as far as we 
    assume that Dawson \cite{dawson} and Takeuti \cite{takeuti1}, \cite{takeuti3} 
    give correct accounts of the development of G\"odel's research in set theory 
    around 1970, we have to  
    say that the description about this episode in \cite{kennedy} contains a couple 
    of inaccuracies:   the note was said to be submitted in 1972 instead of 1970. 
    There is no mention about the continuous effort of G\"odel during  
    1970 to 1975 and instead the article gives the impression that all research 
    of G\"odel was done in the note submitted in 1970. Also Takeuti's 
    contribution is totally suppressed in the article. ]}.
}\fi
Modern treatment of G\"odel's axioms is to be found in Brendle, Larson and 
Todor\v{c}evi\'c \cite{recutangular}. \cite{takeuti2} is cited in  
\cite{recutangular} 
but neither \cite{takeuti1} nor \cite{takeuti3} is mentioned there.

Though the technical details of the present article are rather orthogonal to 
the G\"odel-Takeuti line of the support of $\continuum=\aleph_2$, the first 
author considers the results presented in this article as his belated reply to 
the letter in 1998 and would like to dedicate this article to the 
memory of Professor Gaishi Takeuti. 


\end{document}